\documentclass[11pt]{article}
\usepackage{amsfonts}
\usepackage{amsfonts,mathrsfs,amssymb,amsthm,mathptm}
\usepackage{amsmath,amscd}
\usepackage{mathptm,pslatex}
\usepackage{color}
\usepackage[dvips]{graphicx}
\usepackage[all]{xy}
\usepackage{graphicx}
\usepackage{fancyhdr}

\begin{document}
%\sloppy
\newtheorem{Def}{Definition}[section]
\newtheorem{Bsp}[Def]{Example}
\newtheorem{Prop}[Def]{Proposition}
\newtheorem{Theo}[Def]{Theorem}
\newtheorem{Lem}[Def]{Lemma}
\newtheorem{Koro}[Def]{Corollary}
\theoremstyle{definition}
\newtheorem{Rem}[Def]{Remark}

\newcommand{\add}{{\rm add}}
\newcommand{\gd}{{\rm gl.dim}}
\newcommand{\dm}{{\rm dom.dim}}
\newcommand{\E}{{\rm E}}
\newcommand{\Mor}{{\rm Morph}}
\newcommand{\End}{{\rm End}}
\newcommand{\ind}{{\rm ind}}
\newcommand{\rsd}{{\rm res.dim}}
\newcommand{\rd} {{\rm rep.dim}}
\newcommand{\ol}{\overline}
\newcommand{\overpr}{$\hfill\square$}
\newcommand{\rad}{{\rm rad}}
\newcommand{\soc}{{\rm soc}}
\renewcommand{\top}{{\rm top}}
\newcommand{\pd}{{\rm proj.dim}}
\newcommand{\id}{{\rm inj.dim}}
\newcommand{\fld}{{\rm flat.dim}}
\newcommand{\Fac}{{\rm Fac}}
\newcommand{\Gen}{{\rm Gen}}
\newcommand{\fd} {{\rm fin.dim}}
\newcommand{\DTr}{{\rm DTr}}
\newcommand{\cpx}[1]{#1^{\bullet}}
\newcommand{\D}[1]{{\mathscr D}(#1)}
\newcommand{\Dz}[1]{{\mathscr D}^+(#1)}
\newcommand{\Df}[1]{{\mathscr D}^-(#1)}
\newcommand{\Db}[1]{{\mathscr D}^b(#1)}
\newcommand{\C}[1]{{\mathscr C}(#1)}
\newcommand{\Cz}[1]{{\mathscr C}^+(#1)}
\newcommand{\Cf}[1]{{\mathscr C}^-(#1)}
\newcommand{\Cb}[1]{{\mathscr C}^b(#1)}
\newcommand{\K}[1]{{\mathscr K}(#1)}
\newcommand{\Kz}[1]{{\mathscr K}^+(#1)}
\newcommand{\Kf}[1]{{\mathscr  K}^-(#1)}
\newcommand{\Kb}[1]{{\mathscr K}^b(#1)}
%\stackrel{\sim}
\newcommand{\modcat}{\ensuremath{\mbox{{\rm -mod}}}}
\newcommand{\Modcat}{\ensuremath{\mbox{{\rm -Mod}}}}

\newcommand{\stmodcat}[1]{#1\mbox{{\rm -{\underline{mod}}}}}
\newcommand{\pmodcat}[1]{#1\mbox{{\rm -proj}}}
\newcommand{\imodcat}[1]{#1\mbox{{\rm -inj}}}
\newcommand{\Pmodcat}[1]{#1\mbox{{\rm -Proj}}}
\newcommand{\Imodcat}[1]{#1\mbox{{\rm -Inj}}}
\newcommand{\opp}{^{\rm op}}
\newcommand{\otimesL}{\otimes^{\rm\mathbb L}}
\newcommand{\rHom}{{\rm\mathbb R}{\rm Hom}\,}
\newcommand{\projdim}{\pd}
\newcommand{\Hom}{{\rm Hom}}
\newcommand{\Coker}{{\rm Coker}}
\newcommand{ \Ker  }{{\rm Ker}}
\newcommand{ \Img  }{{\rm Im}}
\newcommand{\Ext}{{\rm Ext}}
\newcommand{\StHom}{{\rm \underline{Hom}}}

\newcommand{\gm}{{\rm _{\Gamma_M}}}
\newcommand{\gmr}{{\rm _{\Gamma_M^R}}}

\def\vez{\varepsilon}\def\bz{\bigoplus}  \def\sz {\oplus}
\def\epa{\xrightarrow} \def\inja{\hookrightarrow}

\newcommand{\lra}{\longrightarrow}
\newcommand{\lraf}[1]{\stackrel{#1}{\lra}}
\newcommand{\ra}{\rightarrow}
\newcommand{\dk}{{\rm dim_{_{k}}}}

\newcommand{\colim}{{\rm colim\, }}
\newcommand{\limt}{{\rm lim\, }}
\newcommand{\Add}{{\rm Add }}
\newcommand{\Tor}{{\rm Tor}}
\newcommand{\Cogen}{{\rm Cogen}}

{\Large \bf
\begin{center}
Stratifications of derived categories from tilting modules
\\
over tame hereditary algebras
%Stratifications of the derived categories of the endomorphism rings \\
%of tilting modules over tame hereditary algebras
\end{center}}
\medskip

\centerline{{\bf Hongxing Chen} and {\bf Changchang Xi$^*$}}
\begin{center} School of Mathematical Sciences, Beijing Normal University, \\
Laboratory of Mathematics and Complex Systems, \\
100875 Beijing, People's Republic of  China \\ E-mail: chx19830818@163.com(H.X.Chen) \quad xicc@bnu.edu.cn(C.C.Xi)\\
\end{center}
\bigskip

\renewcommand{\thefootnote}{\alph{footnote}}
\setcounter{footnote}{-1} \footnote{ $^*$ Corresponding author.
Email: xicc@bnu.edu.cn; Fax: 0086 10 58802136; Tel.: 0086 10
58808877.}
\renewcommand{\thefootnote}{\alph{footnote}}
\setcounter{footnote}{-1} 
\renewcommand{\thefootnote}{\alph{footnote}}
\setcounter{footnote}{-1} 

\begin{abstract} In this paper, we consider the endomorphism algebras
of  infinitely generated tilting modules of the form $R_{\mathcal
U}\oplus R_{\mathcal U}/R$ over tame hereditary $k$-algebras $R$
with $k$ an arbitrary field, where $R_{\mathcal{U}}$ is the
universal localization of $R$ at an arbitrary set $\mathcal{U}$ of
simple regular $R$-modules, and show that the derived module
category of $\End_R(R_{\mathcal U}\oplus R_{\mathcal U}/R)$ is a
recollement of the derived module category $\D{R}$ of $R$ and the
derived module category $\D{{\mathbb A}_{\mathcal{U}}}$ of the
ad\`ele ring ${\mathbb A}_{\mathcal{U}}$ associated with
$\mathcal{U}$. When $k$ is an algebraically closed field, the ring
${\mathbb A}_{\mathcal{U}}$ can be precisely described in terms of
Laurent power series ring $k((x))$ over $k$. Moreover, if $\mathcal
U$ is a union of finitely many cliques, we give two different
stratifications of the derived category of $\End_R(R_{\mathcal
U}\oplus R_{\mathcal U}/R)$ by derived categories of rings, such
that the two stratifications are of different finite lengths.
\end{abstract}

\section{Introduction\label{introduction}}

Tilting modules over tame hereditary algebras have played a special
role in the development of the representation theory of algebras:
Finite-dimensional tilting modules provide a class of minimal
representation-infinite algebras which can be used together with the
covering techniques in \cite{BG} to judge whether an algebra is of
finite representation type or not, while infinite-dimensional
tilting modules involve the generic modules discovered by Ringel in
\cite{R},  Pr\"ufer modules and adic modules. Recently,
Angeleri-H\"ugel and S\'anchez classify all tilting modules over
tame hereditary algebras up to equivalence in \cite{HJ2}. One of the
main ingredients of their classification involves the universal
localizations at simple regular modules, which were already studied
by Crawley-Boevey in \cite{CB}. It is worthy to note that Krause and
Stovicek show very recently in \cite{HO} that over hereditary rings
universal localizations and ring epimorphisms coincide. For
finite-dimensional tilting modules over tame hereditary algebras,
their endomorphism algebras have been well understood from the view
of torsion theory and derived categories (see \cite{BB},
\cite{Happelb}, \cite{hr}, \cite{R1}, and others). Especially, in
this case, there are derived equivalences between the given tame
hereditary algebras and the endomorphism algebras of titling
modules. However, for infinite-dimensional tilting modules, one
cannot get such derived equivalences any more (see \cite{Bz}).
Nevertheless, if they are good tilting modules, then the derived
module categories of their endomorphism rings admit recollements by
derived module categories of the given tame hereditary algebras
themselves on the one side, and of certain universal localizations
of their endomorphism rings on the other side, as shown by a general
result in \cite{CX}. Here, not much is known about the precise
structures of these universal localizations as well as the
composition factors of these recollements. In fact, it seems to be
very difficult to describe them in general.

In the present paper, we will study these new recollements arising
from a class of good tilting modules over tame hereditary algebras
more explicitly. In this special situation, we can describe
precisely the universal localizations appearing in the recollements
in terms of ad\`ele rings which occur often in algebraic number
theory (see \cite[Chapter V]{ne}, determine their derived
composition factors, and provide two completely different
stratifications of the derived module categories of the endomorphism
rings of these tilting modules.

Let $R$ be an indecomposable finite-dimensional tame hereditary
algebra over an arbitrary field $k$. Of our interest are simple
regular $R$-modules. Now, we fix a complete set ${\mathscr S}$ of
all non-isomorphic simple regular $R$-modules, and consider the
equivalence relation $\sim$ on ${\mathscr S}$ generated by
$$L_1\sim L_2 \;\;\mbox{ for}\;\; L_1, L_2\in{\mathscr S}\; \mbox{ if }\;
\Ext^1_R(L_1,L_2)\neq 0.$$ The equivalence classes of this relation
are called cliques (see \cite{CB}). It is well known that all
cliques are finite, and all but at most three cliques consist of
only one simple regular module. For a simple regular $R$-module $L$,
we denote by ${\mathscr C}(L)$ the clique containing $L$. Similarly,
for a subset $\mathcal{V}$ of ${\mathscr S}$, we denote by
${\mathscr C}({\mathcal V})$ the union of all cliques ${\mathscr
C}(L)$ with $L\in\mathcal{V}$.

Let $\mathcal{C}$ be a clique and $V\in\mathcal{C}$. Then there is a
unique Pr\"ufer $R$-module, denoted by $V[\infty]$, such that its
regular socle is equal to $V$ (see \cite{R}).  Moreover, for any two
non-isomorphic simple regular modules in $\mathcal{C}$, the
endomorphism rings of the Pr\"ufer modules corresponding to them are
isomorphic ( see, for instance, Lemma \ref{lem3.3}(3)). Hence we
define $D({\mathcal{C}})$ to be $\End_R(V[\infty])$ for an arbitrary
but fixed module $V\in\mathcal{C}$. It is shown that this ring is a
(not necessarily commutative) discrete valuation ring. Therefore,
the so-called division ring $Q({\mathcal{C}})$ of fractions of
$D({\mathcal{C}})$ exists, which is the ``smallest" division ring
containing $D({\mathcal{C}})$ as a subring up to isomorphism.

Let $\mathcal{U}\subseteq {\mathscr S}$ be a set of simple regular
modules, and let $R_\mathcal{U}$ stand for the universal
localization of $R$ at $\mathcal{U}$ in the sense of Schofield and
Crawley-Boevey. Then it is proved in \cite{HJ} that the $R$-module
$T_{\mathcal U}:= R_{\mathcal{U}}\oplus R_{\mathcal{U}}/R$ is a
tilting module. Following \cite[Example 1.3]{HJ2}, if $\mathcal{U}$
is a union of cliques, the $R$-module $T_{\mathcal{U}}$ is called
the Reiten-Ringel tilting module associated with $\mathcal{U}$. This
class of modules was studied first in \cite{R} and generalized then
in \cite{RR}. As a main objective of the present paper, we will
concentrate us on the derived categories of the endomorphism rings
of tilting modules $T_{\mathcal U}$ for arbitrary subsets $\mathcal
U$ of $\mathscr S$.

Let $k[[x]]$ and $k((x))$ be the algebras of formal and Laurent
power series  over $k$ in one variable $x$, respectively. For an
index set $I$, we define the $I$-ad\`ele ring of $k((x))$ by
$$\mathbb{A}_I:=\bigg\{\big(f_i\big)_{i\in
I}\in\prod_ {i\in I}k((x))\;\big|\;f_i\in k[[x]]\; \mbox{ for almost
all }\; i\in I \bigg\},$$ where $\prod_ {i\in I}k((x))$ stands for
the direct product of $I$ copies of $k((x))$. In particular, if $I$
is a finite set, then ${\mathbb A}_I=k((x))^{|I|}.$

Our main result in this paper is the following theorem, which
provides  us a class of new recollements different  from the one
obtained by the structure of triangular matrix rings.

\begin{Theo}
Let $R$ be an indecomposable finite-dimensional tame hereditary
algebra over an arbitrary field $k$. Let
$\mathcal{U}=\mathcal{U}_0\dot{\cup} \mathcal{U}_1$ be a non-empty
set of simple regular $R$-modules, where $\mathcal{U}_0$ contains no
cliques and $\mathcal{U}_1$ is the union of all cliques $\{{\mathcal
C}_i\}_{i\in I}$ contained in $\mathcal{U}$, where $I$ is an index
set, and let $B$ be the endomorphism ring of $R_{\mathcal{U}}\oplus
R_{\mathcal{U}}/R$, where $R_{\mathcal{U}}$ stands for the universal
localization of $R$ at $\mathcal{U}$. Then there is the following
recollement of derived module categories:

$$
\xymatrix@C=1.2cm{\D{\mathbb{A}_{\mathcal {U}}}\ar[r]&\D{B}\ar[r]
\ar@/^1.2pc/[l]\ar@/_1.2pc/[l]
&\D{R}\ar@/^1.2pc/[l]\ar@/_1.2pc/[l]},\vspace{0.3cm}
$$
where $\mathbb{A}_{\mathcal U}$ is the $I$-ad\`{e}le ring with
respect to the rings $Q({\mathcal{C}_i})$ for $i\in I$, that is,
$$\mathbb{A}_{\,\mathcal{U}}:=\bigg\{\big(f_{i}\big)_{i\in
I}\in\prod_ {i\in I}Q({{\mathcal C}_i})\;\big|\;f_{i}\in
D({{\mathcal C}_i})\; \mbox{\,for almost all }\, i\in I \bigg\}.$$ In
particular, if  $k$ is algebraically closed, then
$\mathbb{A}_{\mathcal U}$ is isomorphic to the $I$-ad\`ele ring
$\mathbb{A}_I$ of the Laurent power series ring $k((x))$.
\label{th1.1}
\end{Theo}

Note that if the field $k$ is algebraically closed then the set
$\mathscr{S}$ of all non-isomorphic simple regular $R$-modules can
be parameterized by the projective line ${\mathbb P}^1(k)$, and the
ad\`ele ring $\mathbb{A}_{\mathbb{P}^1(k)}$ is the same as the
ad\`ele ring $\mathbb{A}_{k(x)}$ of the rational function field
$k(x)$ in global class field theory (see \cite[Chapter VI]{ne} and
\cite[Theorem 2.1.4]{ep} for details). Thus, the ad\`ele ring
$\mathbb{A}_{k(x)}$ occurs in our recollement of Theorem \ref{th1.1}
for the Reiten-Ringel tilting $R$-module $R_{\mathscr{S}}\oplus
R_{\mathscr S}/R$.

As a consequence of Theorem \ref{th1.1}, we can obtain new
stratifications of the derived categories of the endomorphism rings
of tilting modules arising from universal localizations at simple
regular modules.

\begin{Koro}
Let $R$ be an indecomposable finite-dimensional tame hereditary
algebra over an algebraically closed field $k$. Let $r$  be the
number of non-isomorphic simple $R$-modules. Suppose that
$\mathcal{U}$ is a non-empty finite subset of ${\mathscr S}$
consisting of $s$ cliques. Let $B$ be the endomorphism ring of the
Reiten-Ringel tilting $R$-module associated with ${\mathcal{U}}$.
Then $\D{B}$ admits two stratifications by derived module
categories, one is of length $r+s$ with the composition factors: $r$
copies of the ring $k$ and $s$ copies of the ring $k((x))$, and the
other is of length $r+s-1$ with the composition factors: $r-2$
copies of the ring $k$, $s$ copies of the ring $k[[x]]$ and one copy
of a Dedekind integral domain contained in the ring $k(x)$.
\label{cor1.2}
\end{Koro}

Observe that if $R$ is the Kronecker algebra and $\mathcal{U}$
consists of one simple regular module, then we re-obtain the
stratifications, shown in the example of \cite[Section 8]{CX}, from
Corollary \ref{cor1.2}.

\medskip
Now, let us state the structure of this paper. In Section 2, we fix
notations  and recall some definitions and basic facts which will be
used throughout the paper. In Section 3, we first prepare with a few
lemmas, and then prove the main result, Theorem \ref{th1.1}. In
Section 4, we first consider the  case of general tame hereditary
algebras, and then turn to the special case of the Kronecker
algebra. With these preparations in hand, together with a result in
\cite{hx2}, we can determine the derived composition factors of the
derived categories of the endomorphism rings of Reiten-Ringel
tilting modules, and therefore get a proof of Corollary
\ref{cor1.2}.

\medskip
{\bf Acknowledgements}. The author H. X. Chen would like to thank
the Doctor Funds of Beijing Normal University for partial support,
and the corresponding author C. C. Xi would like to acknowledge
partial support by PCSIRT. The paper is revised during a visit of Xi
to the Hausdorff Research Institute for Mathematics, Bonn, Germany
in February and March, 2011 for a trimester program, he would like
to thank C. M. Ringel for discussion on the subject, and the
organizers for invitation.

\section{Preliminaries\label{sect2}}

First, we recall some standard notations which will be used
throughout this paper.

All rings considered are assumed to be associative and with
identity, all ring homomorphisms preserve identity, and all full
subcategories $\mathcal D$ of a given category $\mathcal C$ are
closed under isomorphic images, that is, if $X$ and $Y$ are objects
in $\cal C$, then $Y\in{\mathcal D}$ whenever $Y\simeq X$ with $X\in
{\mathcal D}$.

Let $R$ be a ring.

We denote by $R\Modcat$ the category of all unitary left
$R$-modules, and by $R\modcat$ the category of finitely generated
unitary left $R$-modules. Unless stated otherwise, by an $R$-module
we mean a left $R$-module. For an $R$-module $M$, we denote by
$\add(M)$ \big(respectively, $\Add(M)$\big) the full subcategory of
$R\Modcat$ consisting of all direct summands of finite
(respectively, arbitrary) direct sums of copies of $M$. If $I$ is an
index set, we denote by $M^{(I)}$ the direct sum of $I$ copies of
$M$.

If $f: M\ra N$ is a homomorphism of $R$-modules, then the image of
$x\in M$ under $f$ is denoted by $(x)f$ instead of $f(x)$. Also, for
any $R$-module $X$, the induced morphisms $\Hom_R(X,f):
\Hom_R(X,M)\ra \Hom_R(X,N)$ and $\Hom_R(f,X): \Hom_R(N, X)\ra
\Hom_R(M, X)$ are denoted by $f^*$ and $f_*$, respectively.

Given a class $\mathcal U$ of $R$-modules, we denote by ${\mathcal
F}({\mathcal U})$ the full subcategory of $R$-Mod consisting of all
those $R$-modules $M$ which have a finite filtration $0=M_0\subseteq
M_1\subseteq \cdots \subseteq M_n=M$ such that $M_{i}/M_{i-1}$ is
isomorphic to a module in $\mathcal U$ for each $i$. We say that $M$
is a direct union of finite extensions of modules in $\mathcal U$ if
$M$ is the direct limit of a direct system of submodules of $M$
belonging to ${\mathcal F}({\mathcal U})$.

Let $\D{R}$ be the (unbounded) derived category of $R\Modcat$, which
is the localization of the homotopy category of $R\Modcat$ at all
quasi-isomorphisms. Furthermore, we always identify $R\Modcat$ with
the full subcategory of $\D{R}$ consisting of all stalk complexes
concentrated on degree zero. It is well known that
$\Hom_{\D{R}}(X,Y[n])\simeq\Ext^n_R(X,Y)$ for any $X,Y\in R\Modcat$
and $n\in \mathbb{N}$, where $[n]$ stands for the $n$-th shift
functor of $\D{R}$, and that the triangulated category $\D{R}$ has
small coproducts, that is, coproducts indexed by sets exist in
$\D{R}$.

If $R$ is  an Artin $k$-algebra over a commutative Artin ring $k$,
we denote by $D$ the usual duality, and by $\tau$ the
Auslander-Reiten translation of $R$.

Now, let us recall the definition of recollements of triangulated
categories. This notion was first introduced by Beilinson, Bernstein
and Deligne in \cite{BBD} to study the triangulated categories of
perverse sheaves over singular spaces, and later was used by Cline,
Parshall and Scott in \cite{CPS} to stratify the derived categories
of quasi-hereditary algebras arising from the representation theory
of semisimple Lie algebras and algebraic groups.

Let $\mathcal{D}$ be a triangulated category with small coproducts.
We denote by $[1]$ the shift functor of  $\mathcal{D}$.

\begin{Def}\rm \cite{BBD}
Let $\mathcal{D'}$ and $\mathcal{D''}$ be triangulated categories.
We say that $\mathcal{D}$ is a recollement of $\mathcal{D'}$ and
$\mathcal{D''}$ if there are six triangle functors $i_*, i^*,i^!,
j^!, j_*$ and $j_!$ as in the following diagram
$$\xymatrix{\mathcal{D''}\ar^-{i_*=i_!}[r]&\mathcal{D}\ar^-{j^!=j^*}[r]
\ar_-{i^!}@/^1.8pc/[l]\ar_-{i^*}@/_1.8pc/[l]
&\mathcal{D'}\ar_-{j_*}@/^1.8pc/[l]\ar_-{j_!}@/_1.8pc/[l]}$$ such
that

$(1)$ $(i^*,i_*),(i_!,i^!),(j_!,j^!)$ and $(j^*,j_*)$ are adjoint
pairs,

$(2)$ $i_*,j_*$ and $j_!$ are fully faithful,

$(3)$ $i^!j_*=0$ (and thus also $j^! i_!=0$ and $i^*j_!=0$),

$(4)$ for each object $C\in\mathcal{D}$, there are two triangles in
$\mathcal D$:
$$
i_!i^!(C)\lra C\lra j_*j^*(C)\lra i_!i^!(C)[1] \quad \mbox{and}
\quad j_!j^!(C)\lra C\lra i_*i^*(C)\lra j_!j^!(C)[1].
$$
\label{def2.1}
\end{Def}

In the following, if $\mathcal{D}$ is a recollement of
$\mathcal{D'}$ and $\mathcal{D''}$, we also say that there is a
recollement among $\mathcal{D'}$, $\mathcal{D}$ and $\mathcal{D''}$,
or very briefly, that $\mathcal{D}$ admits a recollement.

A well known example of recollements of derived categories of rings
is given by triangular matrix rings: Suppose that $A$, $B$ are
rings, and that
$M$ is an $A$-$B$-bimodule. Let  $R=\left(\begin{array}{lc} A & M\\
0 & B\end{array}\right)$ be the triangular matrix ring associated
with $A,B$ and $M$. Then there is a recollement of derived
categories:

$$\xymatrix@C=1.2cm{\D{A}\ar[r]&\D{R}\ar[r]
\ar@/^1.2pc/[l]\ar@/_1.2pc/[l]
&\D{B}\ar@/^1.2pc/[l]\ar@/_1.2pc/[l]}.\vspace{0.3cm}$$ A
generalization of this situation is the so-called stratifying ideals
defined by Cline, Parshall and Scott, and can be found in
\cite{CPS}.

Another type of examples of recollements of derived categories of
rings appears in the tilting theory of infinitely generated tilting
modules over arbitrary rings (see \cite{CX}). Before we state this
kind of examples, we recall first the definition of tilting modules
over arbitrary rings from \cite{ct}, and then the notion of
universal localizations which is closely related to constructing
tilting modules.

\begin{Def}\rm \label{def2.2}
An  $R$-module $T$ is called a tilting module (of projective
dimension at most one) if the following conditions are satisfied:

$(T1)$ The projective dimension of $T$ is at most $1$, that is,
there exists an exact sequence: $0\rightarrow P_1\rightarrow
P_0\rightarrow T\rightarrow 0$ with $P_i$ projective for $i=0,1$,

$(T2)$ $\Ext^i_R(T,T^{({\alpha})})=0$ for each $i\geq 1$ and each
index set $\alpha$,  and

$(T3)$ there exists  an exact sequence
$$0\lra {}_RR\lra
T_0\longrightarrow T_1\longrightarrow 0$$ of $R$-modules such that
$T_i\in\Add(T)$ for $ i=0,1$.

A tilting $R$-module $T$ is called good if $T_0$ and $ T_1$ in (T3)
lie in $\add(T)$, and classical if $T$ is good and finitely
presented.
\end{Def}

A special kind of good tilting modules can be constructed from
injective ring epimorphisms, including particularly certain
universal localizations. The following result on universal
localizations is well known.

\begin{Lem}{\rm (see \cite{cohenbook1}, \cite{Sch1})} Let $R$ be a ring and
$\Sigma$ a set of homomorphisms between finitely generated
projective $R$-modules. Then there is a ring $R_{\Sigma}$ and a
homomorphism  $\lambda: R\to R_{\Sigma}$ of rings  with the
following properties:

$(1)$ $\lambda$ is $\Sigma$-inverting, that is, if $\alpha:P\to Q$
belongs to $\Sigma$, then
$R_{\Sigma}\otimes_R\alpha:R_{\Sigma}\otimes_RP\to
R_{\Sigma}\otimes_R Q$ is an isomorphism of $R_{\Sigma}$-modules,
and

$(2)$ $\lambda$ is universal $\Sigma$-inverting, that is, if $S$ is
a ring such that there exists a $\Sigma$-inverting homomorphism
$\varphi:R\to S$, then there exists a unique homomorphism
$\psi:R_{\Sigma}\to S$ of rings such that $\varphi=\lambda\psi$.

$(3)$ The homomorphism $\lambda:R\to R_{\Sigma}$  is a ring
epimorphism with $\Tor^R_1(R_{\Sigma}, R_{\Sigma})=0.$
\label{lem2.3}
\end{Lem}

We call  $\lambda: R\to R_{\Sigma}$ in Lemma \ref{lem2.3} the
universal localization of $R$ at $\Sigma$. Recall that, by
\cite[Theorem 2.5]{HJ}, if $\lambda$ is injective and the $R$-module
$R_\Sigma$ has projective dimension at most one, then
$R_\Sigma\oplus R_\Sigma/R$ is a tilting $R$-module.

Of particular interest are the following two kinds of universal
localizations.

The first one is associated with subsets of elements in rings. Let
$\Phi$ be a non-empty subset of $R$. Then we consider the universal
localization of $R$ at all homomorphisms $\rho_r$ with $r\in\Phi$,
where $\rho_r$ is the right multiplication map $R\to R$ defined by
$x\mapsto xr$ for $x\in R$. For simplicity, we write $R_{\Phi}$ for
this universal localization, and say that $R_\Phi$ is the universal
localization of $R$ at $\Phi$. Note that, by the property of
universal localizations, $R_\Phi$ is also isomorphic to the ``right"
universal localization of $R$ at all left multiplication maps
$\sigma_r: R_R\to R_R$ defined by $x\mapsto rx$ for $x\in\Phi$,
which are regarded as homomorphisms of right $R$-modules. Clearly,
if $0\in\Phi$, then $R_\Phi=0$. If $0\notin\Phi$, then we consider
the smallest multiplicative subset of $R$ containing $\Phi$, and get
$R_\Phi=R_{\Phi_1}$. Recall that a subset $\Phi$ of $R$ is said to
be multiplicative if $0\notin \Phi$, $1\in\Phi$, and it is closed
under multiplication.

From now on, we assume that $\Phi$ is a multiplicative subset of
$R$.

Under some extra assumptions on $\Phi$, the ring $R_\Phi$ can be
characterized  by Ore localizations which generalizes the notion of
localizations in commutative rings. To explain this point in detail,
we first recall some relevant definitions about Ore localizations.
For more details, we refer to \cite[Chapter 4]{Lam}.

\begin{Def}\label{denom} {\rm
A subset $\Phi$ of $R$ is called a left denominator subset of $R$ if
$\Phi$ satisfies the following two conditions: $(i)$ For any $a\in
R$ and $s\in\Phi$, there holds $\Phi a\cap Rs\neq \emptyset$, and
$(ii)$ for any $r\in R$, if $rt=0$ for some $t\in\Phi$, then there
exists some $t'\in\Phi$  such that $t'r=0$. If $\Phi$ satisfies only
the condition $(i)$, then $\Phi$ is called a left Ore subset of
$R$.}
\end{Def}

Similarly, we can define the notions of right denominator sets and
right Ore sets, respectively. Clearly, if $R$ is commutative, then
every multiplicative subset of $R$ is  a left and right denominator
set. Furthermore, if $R$ is a domain, that is, $R$ is a (not
necessarily commutative) ring which has neither left zero-divisors
nor right zero-divisors, then $R\setminus \{0\}$ is a left
denominator set if and only if it is a left Ore set if and only if
$Rr_1\cap Rr_2\neq \{0\}$ for any non-zero elements $r_1,r_2\in R$.
We say that $R$ is a left Ore domain if $R\setminus \{0\}$ is a left
denominator set.

The following lemma explains how left Ore localizations arise, and
establishes a relationship between left Ore localizations and
universal localizations.

\begin{Lem}{\rm \cite[Theorem 10.6, Corollary
10.11]{Lam}} Let $\Phi$ be a left denominator subset of $R$ and
$\lambda: R\to R_\Phi$ the universal localization of $R$ at $\Phi$.
Then there is a ring, denoted by $\Phi^{-1}R$, and a ring
homomorphism $\mu:R\to \Phi^{-1}R$ such that

$(1)$ $\mu$ is $\Phi$-invertible, that is, $(s)\mu$ is a unit in
$\Phi^{-1}R$ for each $s\in\Phi$,

$(2)$ every element of $\,\Phi^{-1}R$ has the form
$\big((t)\mu\big)^{-1}(r)\mu$ for some $t\in\Phi$ and some $r\in R$,

$(3)$ $\ker(\mu)=\{r\in R\mid sr=0\; \mbox{for\, some}\;
s\in\Phi\}$, and

$(4)$ there is a unique isomorphism $\nu: \Phi^{-1}R\to R_\Phi$ of
rings such that $\lambda=\mu\nu$. \label{Ore}
\end{Lem}

The ring $\Phi^{-1}R$ in Lemma \ref{Ore} is called a left ring of
fractions of $R$ (with respect to $\Phi\subseteq R$), or
alternatively, a left Ore localization of $R$ at $\Phi$. Clearly,
for commutative rings, Ore localizations  and the usual
localizations at multiplicative subsets coincide.

Similarly, when $\Phi$ is a right denominator subset of $R$, we can
define a right ring $R\Phi^{-1}$ of fractions of $R$. If $\Phi$ is a
left and right denominator subset of $R$, then $\Phi^{-1}R$ is
called the ring of fractions of $R$, or the Ore localization of $R$
at $\Phi$. Actually, in this case, both $\Phi^{-1}R$ and
$R\Phi^{-1}$ are isomorphic to $R_\Phi$. Furthermore, if $R$ is a
left and right Ore domain $R$, then the ring of fractions of $R$
with respect to $R\setminus \{0\}$ is usually denoted by $Q(R)$.
Notice that, up to isomorphism, $Q(R)$ is the smallest division ring
containing $R$ as a subring. So we call $Q(R)$ the division ring of
fractions of $R$.

The other kind of universal localizations is provided by universal
localizations at injective homomorphisms between finitely generated
projective modules, and therefore  related to finitely presented
modules of projective dimension at most one.

Suppose that $\mathcal {U}$ is a set of finitely presented
$R$-modules of projective dimension at most one. For each
$U\in\mathcal{U}$, there is a finitely generated projective
presentation of $U$, that is, an exact sequence of $R$-modules
$$
(\ast)\quad 0\lra P_{1}\lraf{f_U} P_0\lra U\lra 0,
$$
such that $P_1$ and $P_0$ are finitely generated and projective. Set
$\Sigma:=\{f_U\mid U\in\mathcal{U}\}$, and let $R_\mathcal{U}$ be
the universal localization of $R$ at $\Sigma$.  If $f'_U: Q_1\to
Q_0$ is another finitely generated projective presentation of $U$,
then the universal localization of $R$ at $\Sigma':=\{f'_U\mid U\in
{\mathcal U}\}$ is isomorphic to $R_\mathcal{U}$. Hence
$R_\mathcal{U}$ does not depend on the choices of the injective
homomorphisms $f_U$, and we may say that $R_{\mathcal U}$ is the
universal localization of $R$ at $\mathcal U$.

Clearly, we have $\Tor_i^R(R_{\mathcal U},U)=0$ for all $i\ge 0$ and
$U\in \mathcal{U}$, and therefore $\Tor_i^R(R_{\mathcal U},X)=0$ for
all $i\ge 0$ and $X\in \mathcal{F}(\mathcal{U})$.

\medskip
Now, we state the promised example of recollements as a proposition
which is a consequence of \cite[Lemma 6.2, Corollary 6.6]{CX}. It is
worthy to notice that the recollement in this proposition is, in
general, different from the one obtained from the structure of
triangular matrix rings.

\begin{Prop} Let $\mathcal {U}$ be a set of finitely presented $R$-modules of
projective dimension one, and let $\lambda: R\to R_\mathcal{U}$ be
the universal localization  of $R$ at $\mathcal{U}$. Suppose that
$\lambda$ is injective and that the $R$-module $R_\mathcal{U}$ has
projective dimension at most one. Set $S:=\End_R(R_\mathcal{U}/R)$,
$B:=\End_R(R_\mathcal{U}\oplus R_\mathcal{U}/R)$ and
$\Sigma:=\{S\otimes_R f_U\mid U\in\mathcal{U}\}$. Then there is a
recollement of derived module categories:

$$\xymatrix@C=1.2cm{\D{S_{\Sigma}}\ar[r]&\D{B}\ar[r]
\ar@/^1.2pc/[l]\ar@/_1.2pc/[l]
&\D{R}\ar@/^1.2pc/[l]\ar@/_1.2pc/[l]},\vspace{0.3cm}$$ where
$S_\Sigma$ is the universal localization  of $S$ at $\Sigma$.
\label{lem2.4}
\end{Prop}

In many cases we can use this proposition repeatedly because the
following result states that iterated universal localizations are
again universal localizations.

\begin{Lem} {\rm \cite[Theorem 4.6]{Sch1}}
Let $\Sigma$ and $\Gamma$ be sets of homomorphisms between finitely
generated projective $R$-modules. Set
$\ol{\Gamma}:=\{R_{\Sigma}\otimes_R f\mid f\in\Gamma\}$. Then the
universal localization of $R$ at $\Sigma\cup\Gamma$ is isomorphic to
the universal localization of $R_\Sigma$ at $\ol{\Gamma}$, that is,
$R_{\Sigma\cup\Gamma}\simeq (R_{\Sigma})_{\ol{\Gamma}}$ as rings.
\label{lem2.5}
\end{Lem}

Next, we recall the definition of discrete valuation rings.

\begin{Def} \rm \label{def2.3}
A ring $R$ is called a  discrete valuation ring (which may not be
commutative) if the following conditions hold true:

$(1)$ $R$ is a local ring, that is, $R$ has a unique maximal left
ideal $\mathfrak m$;

$(2)$ $\bigcap_{i\geq 1}\mathfrak {m}^i=0$;

$(3)$ $\mathfrak {m}=pR=Rp$, where $p$ is some non-nilpotent element
of $R$.
\end{Def}

We remark that an equivalent definition of discrete valuation rings
is the following: A non-division ring $R$ is called a discrete
valuation ring if it is a local domain with $\mathfrak{m}$ the
unique maximal ideal of $R$ such that the only left ideals and the
only right ideals of $R$ are of the form $\mathfrak{m}^i$ for
$i\in\mathbb{N}$.

The element $p$ in the above condition $(3)$ is called a prime element of
$R$. Clearly, for each invertible element $v$ of $ R$,  both $vp$
and $pv$ are prime elements. A discrete valuation ring is said to be
complete if the canonical map
$R\to\varprojlim_i\,R/{\mathfrak{m}^i}$ is an isomorphism. Note that
every discrete valuation ring can be embedded into a complete
discrete valuation ring.

The following lemma collects some basic properties of discrete
valuation rings, which will be frequently used in our proofs.

\begin{Lem}{\rm (\cite[Chapter 1]{KT}, \cite{Lam})}\label{dv}
Let $R$ be a discrete valuation ring, $\mathfrak m$ the unique
maximal ideal of $R$, and $p$ a prime element of $R$. Then the
following statements are true:

$(1)$ The ideals $\mathfrak{m}^i$ $(i\in\mathbb{N})$ are the only
left ideals and the only right ideals of $R$.

$(2)$ For any non-zero element $x\in R$, there are unique elements
$x_1, x_2\in R\setminus\mathfrak{m}$ such that $x=x_1p^n=p^nx_2$ for
some $n\in\mathbb{N}$.

$(3)$  $R$ is a left and right Ore domain. In particular, the
division ring $Q(R)$ of fractions of $R$ exists.

$(4)$ $Q(R)$ is isomorphic to the universal localization of $R$ at
 the map $\rho_p: R\to R$ defined by $r\mapsto rp$ for $r\in R$.
\end{Lem}

\medskip
Finally, we prepare several homological results for our later
proofs.

\begin{Lem}
Let $R$ be a ring  and let $0\lra X \lraf{(f,\,g)} Y\oplus Z\lraf{h}
W\lra 0$ be an exact sequence of $R$-modules. Assume that $f: X\ra
Y$ is injective and that there is a homomorphism $\tilde{g}: Y\ra Z$
with $g=f\tilde{g}: X\ra Z$. Then there exists an automorphism
$\varphi$ of the module $Y\oplus Z$ and an isomorphism $\psi:W\ra
\Coker(f)\oplus Z$ such that the following diagram commutes:
$$\xymatrix{0\ar[r]&X \ar[r]^-{(f,\,g)}\ar@{=}[d]& Y\oplus
Z\ar[r]^-{h}\ar[d]^-{\varphi}& W\ar[r]\ar[d]^-{\psi}& 0\\
0\ar[r]& X\ar[r]^-{(f,\,0)}& Y\oplus Z\ar[r]^-{\left({\pi\;
0}\atop{0\;1}\right)}& \Coker(f)\oplus Z\ar[r]& 0, }$$ where $\pi:
Y\ra\Coker(f)$ stands for the canonical surjection. \label{lem3.6}
\end{Lem}

{\it Proof.} Set $\varphi:={\left({\;\;\,1\quad
{-\tilde{g}}}\atop{0\quad 1}\right)}$. Then $\varphi$ is an
automorphism of the module $Y\oplus Z$. Since $g=f\tilde{g}$, we
have $(f, g)\varphi=(f,0)$. Thus, there exists a unique homomorphism
$\psi:W\lra \Coker(f)\oplus Z$, such that the exact diagram in Lemma
\ref{lem3.6} is commutative. Clearly, $\psi$ is an isomorphism. This
completes the proof. $\square$

\medskip
The following homological facts are well known in the literature
(see, for example, the book \cite{gt}).

\begin{Lem} Let $R$ be a ring.

$(1)$ If $\{X_{\alpha}\}_{\alpha\in I}$ is a direct system of
$R$-modules, then

\smallskip
(i) $\Hom_R(\displaystyle\varinjlim_{\alpha} X_{\alpha},M)\simeq
\displaystyle\varprojlim_{\alpha} \Hom_R(X_{\alpha},M)$ for any
$R$-module $M$.

\smallskip
(ii) For any finitely presented $R$-module $M$, we have
$\Hom_R(M,\displaystyle\varinjlim_{\alpha}
X_{\alpha})\simeq\displaystyle\varinjlim_{\alpha}\Hom_R(M,X_{\alpha})$.

(iii) Let $n\ge 0$. If $M$ is an $R$-module with  a projective
resolution $ \cdots \ra P_{n+1}\ra \cdots\ra P_1\ra P_0\ra M\ra 0$
such that all $P_j$, with $0\le j\le n+1$, are finitely generated,
then
$$\Ext^i_R(M,\displaystyle\varinjlim_{\alpha} X_{\alpha})\simeq\displaystyle\varinjlim_{\alpha}\Ext^i_R(M,X_{\alpha})$$
for all $i\le n$.

(iv) If $M$ is a pure-injective $R$-module (for example, $M$ is of
finite length over its endomorphism ring), then
$$\Ext^i_R(\displaystyle\varinjlim_{\alpha} X_{\alpha},M)\simeq
\displaystyle\varprojlim_{\alpha}\Ext^i_R(X_{\alpha}, M)$$
for all $i\ge 0$. Conversely, if this isomorphism is true for $i=1$
and for every directed system $X_{\alpha}$, then $M$ is
pure-injective.

$(2)$ If $\{Y_{\alpha}\}_{\alpha\in I}$ is an inverse system of
$R$-modules, then, for any $R$-module $M$,
$$\Hom_R(M,\displaystyle\varprojlim_{\alpha} Y_{\alpha})\simeq
\displaystyle\varprojlim_{\alpha}\Hom_R(M, Y_{\alpha}).$$
\label{2.7}
\end{Lem}

{\it Remarks}. (1) The statement (iv) is due to Maurice Auslander.

(2) The class of all pure-injective $R$-modules is closed under
products, direct summands and finite direct sums. In general, it is
not closed under extensions.

\begin{Lem} Let $A$ be a finite-dimensional $k$-algebra over a field $k$,
$M$  a finite-dimensional $A$-module and
$N$ an arbitrary $A$-module.

(1) If $\pd(M)\le 1$, then $D\Ext^1_A(M,N)\simeq \Hom_A(N,\tau M)$,
where $\pd(M)$ stands for the projective dimension of $M$.

(2) If $\id(M)\le 1$, then $\Ext^1_A(N,M)\simeq D\Hom_A(\tau^{-1}M,
N))$, where $\id(M)$ stands for the injective dimension of $M$.
\label{2.8}
\end{Lem}

{\it Proof.} It is known that every $A$-module $N$ is a direct limit
of finitely presented $A$-modules $\{X_{\alpha}\}_{\alpha\in I}$,
and that (1) and (2) hold true for finitely generated modules $N$.
Then, it follows from Lemma \ref{2.7} that
$$\begin{array}{rl}
D\Ext^1_A(M,N) & \simeq D\Ext^1_A(M,\displaystyle\varinjlim_{\alpha}
X_{\alpha})\simeq
D\,\displaystyle\varinjlim_{\alpha}\Ext^1_A(M,X_{\alpha}) \simeq\displaystyle\varprojlim_{\alpha}D\Ext^1_A(M,X_{\alpha})\\
& \simeq \displaystyle\varprojlim_{\alpha}\Hom_A(X_{\alpha}, \tau M)
\simeq \Hom_A(\displaystyle\varinjlim_{\alpha} X_{\alpha},\tau M)=
\Hom_A(N,\tau M). \end{array}$$ This proves (1).  The statement (2)
can be shown similarly. $\square$

\section{Proof of the main result \label{sect3}}
Unless stated otherwise, we assume  from now on that $R$ is  an
indecomposable finite-dimensional tame hereditary algebra over an
arbitrary but fixed field $k$.

Let ${\mathscr S}:={\mathscr S}\!(R)$ be a fixed complete set of
isomorphism classes of all simple regular $R$-modules. For each
$U\in{\mathscr S}$ and $n>0$, we denote by $U[n]$ the $R$-module of
regular length $n$ on the ray
$$(\ast)\quad U=U[1]\subset U[2]\subset\cdots\subset U[n]\subset
U[n+1]\subset\cdots,$$ and let $U[\infty]=
\displaystyle\varinjlim_n\; U[n]$ be the Pr\"ufer module
corresponding to $U$. Note that $U[\infty]$ has a unique regular
submodule $U[n]$ of regular length $n$, and therefore admits a
unique chain of regular submodules, and that each endomorphism of
$U[\infty]$ restricts to an endomorphism of $U[n]$ for any $n>0$.
For further information on regular modules and Pr\"ufer modules over
tame hereditary algebras, we refer to \cite[Section 4, 5]{R} and
\cite{dr}.

Recall that we have defined an equivalence relation $\sim$ on
$\mathscr{S}$ in Section \ref{introduction}.  It is known that two
simple regular modules lie in the same clique if and only if they
lie in the same tube. Thus a clique is just the set of all simple
regular modules belonging to a fixed tube.

Let $U\in{\mathscr S}$ and  $\mathcal{U}\subseteq {\mathscr S}$. We
denote by ${\mathscr C}(U)$ the clique containing $U$, and by $c(U)$
the cardinality of ${\mathscr C}(U)$. Similarly, we denote by
${\mathscr C}({\mathcal{U}})$ the union of all cliques ${\mathscr
C}(U)$ with $U\in\mathcal{U}$, and by $c({\mathcal U})$ the
cardinality of ${\mathscr C}({\mathcal U})$. As mentioned before,
$c(U)$ is always finite, and furthermore, $c(U)=1$ for almost all
$U\in{\mathscr S}$. In fact, there are at most $3$ cliques
consisting of more than one element. Also, we know that all cliques
consist of one simple regular $R$-module if and only if $R$ has only
two isomorphism classes of simple modules. If $k$ is an
algebraically closed field, this is  equivalent to that $R$ is
Morita equivalent to the Kronecker algebra.

\subsection{Endomorphism rings of direct sums of Pr\"ufer modules}

In this subsection, we shall consider the endomorphism ring of the
direct sum of all Pr\"ufer modules obtained from a given tube. This
ring was calculated already in \cite{R}. For convenience of the
reader and also for  the later proof of our main result, we include
here some details of this calculation.

Throughout this subsection, let $\mathcal{C}$ be a clique of
$R$-mod, $U\in\mathcal{C}$, and $\bf{t}$ the tube of rank $m\geq1$
containing $\mathcal{C}$. Set $U_i:=\tau^{-(i-1)}U$ for
$i\in\mathbb{Z}$. Then $\tau^{-m}U\simeq U$, and
$\mathcal{C}=\{U_1,U_2,\cdots, U_{m-1}, U_m\}$ which is a complete
set of non-isomorphic simple regular modules in $\bf{t}$. Since
$U_j\simeq U_{j+m}$ for any $j\in\mathbb{Z}$, the subscript of $U_j$
is always modulo $m$ in our discussion below. It is well known that
$\End_R(U_i)$ is a division algebra and $\Hom_R(U_i,U_j)=0$ for
$1\le i\neq j\le m$, and that $D\Ext^1_R(U_i,U_j)\simeq \End_R(U_i)$
if $j=i-1$, and zero otherwise. Furthermore, $\bf{t}$ is an exact
abelian subcategory of $R\modcat$, and every indecomposable module
in $\bf{t}$ is serial, that is, it has a unique regular composition
series in $\bf{t}$. For example, for any $i\in\mathbb{Z}$ and $j>0$,
the module $U_i[j]$ admits successive regular composition factors
$U_i,U_{i+1},\cdots,U_{i+j-1}$ with $U_i$ as its unique regular
socle and with $U_{i+j-1}$ as its unique regular top. For details,
see \cite[Section 3.1]{R1}.

Now, we mention some properties of Pr\"ufer modules.

\begin{Lem} The following statements hold true for the tube $\bf{t}$.

$(1)$ For any $1\le i\le m$ and for any regular module $X$ in ${\bf
t}$, we have $\Hom_R(U_i[\infty], X)=0=\Ext^1_R(X,U_i[\infty])$.
Further, if $1\le i<j\le m$, then $\Hom_R(U_i[n],U_j[\infty])=0$ for
$1\le n \le j-i$, and $\Hom_R(U_j[n],U_i[\infty])=0$ for $1\le n\le
m-j+i$.

$(2)$ Let $i,j\in{\mathbb N}$ with $1\leq i<j$. Then, for any
$n>j-i$, there is a canonical exact sequence of $R$-modules:
$$0\lra U_i[j-i]\lra U_i[n]\lraf{\epsilon_{i,j}[n]} U_{j}[n-(j-i)]\lra 0.$$
In particular, we get a canonical exact sequence
$$0\lra U_i[j-i]\lra U_i[\infty]\lraf{\epsilon_{i,j}}U_{j}[\infty]\lra 0, $$
where $\epsilon_{i,j}:=
\displaystyle\varinjlim_n\,\epsilon_{i,j}[n]$. Moreover,  we have
$\epsilon_{i,j}=\epsilon_{i+m,\,j+m}$  and
$\epsilon_{i,j}\,\epsilon_{j,p}=\epsilon_{i,p}$ for any $p>j$.

 $(3)$ Let $i,j\in{\mathbb N}$ with $1\leq j-i<m$. Then $\epsilon_{i,j}$
induces an isomorphism of left $\End_R(U_i[\infty])$-modules:
$$(\epsilon_{i,j})^*:\End_R(U_i[\infty])\lraf{\sim}\Hom_R(U_i[\infty],U_j[\infty]),$$
and an isomorphism of right $\End_R(U_j[\infty])$-modules:
$$(\varepsilon_{i,j})_*:\End_R(U_j[\infty])\lraf{\sim}\Hom_R(U_i[\infty],U_j[\infty]).$$
In  particular, we get a ring isomorphism
$\varphi_{i,j}:\End_R(U_i[\infty])\ra \End_R(U_j[\infty])$, $
f\mapsto f'$ for $f\in \End_R(U_i[\infty])$ and $f'\in
\End_R(U_j[\infty])$,with $f\epsilon_{i,j}=\epsilon_{i,j}f'$.

$(4)$ Suppose $1\leq r,\,s,\, t\leq m$. Set $\Delta_{r,\,s}:=
\left\{\begin{array}{ll} 0 & \mbox{if } r<s,\\
1 & \mbox{if }r\geq s,\end{array} \right.$ and define $\pi_{r,
s}:=\epsilon_{r,\,s+\Delta_{r,s}\,m}\in \Hom_R(U_r[\infty],
U_{s+\Delta_{r,s}m}\,[\infty])$. Then
$$\pi_{r,s}\,\pi_{s,t}=
\left\{\begin{array}{ll} \pi_{r,t} & \mbox{if } \Delta_{r,s}+\Delta_{s,t}=\Delta_{r,t},\\
\pi_{r,r}\,\pi_{r,t} &\mbox{otherwise}.\end{array} \right.$$ In
particular, we have $(\pi_{i, i})\varphi_{i,j}=\pi_{j,j}$ for any
$1\leq i<j\leq m$.

$(5)$ The ring $\End_R(U_i[\infty])$ is a  complete discrete
valuation ring with $\pi_{i,i}$ as a prime element. If $k$ is an
algebraically closed field,  then there is a ring isomorphism
$\varphi_i:\End_R(U_i[\infty])\to k[[x]]$ which sends $\pi_{i,i}$ to
$x$.
 \label{lem3.3}
\end{Lem}

{\it Proof.}  $(1)$ Note that $D\Ext^1_R(X,U_i[\infty])\simeq
\Hom_R(U_i[\infty],\tau X)$ for any $X\in {\bf t}$ by Lemma
\ref{2.8}(1), and that every indecomposable module in $\bf{t}$ is
serial. This means that, to prove the first statement in (1), it
suffices to show $\Hom_R(U_i[\infty],U_j)=0$ for all $1\le j\le m$.
In fact, since the inclusion map $U_i[n]\ra U_i[n+1]$ induces a zero
map from $\Hom_R(U_i[n+1],U_j)$ to $ \Hom_R(U_i[n],U_j)$ for all
$n$. This implies that
$$\Hom_R(U_i[\infty],U_j)= \Hom_R(\displaystyle\varinjlim_n
U_i[n],U_j)\simeq \displaystyle\varprojlim_n\,\Hom_R(U_i[n],U_j)=
0.$$
%Observe that, for another different proof of $\Ext^1_R(X,U_i[\infty])=0$, we refer to \cite[Section 4.5]{R}.

The last statement in (1) follows from the fact that the abelian
category $\bf t$ is serial.

$(2)$ For any $n>j-i$, we can easily see from the structure of the
tube $\bf{t}$ that there is an  exact commutative diagram of
$R$-modules:
$$\xymatrix{
0\ar[r]&U_i[j-i]\ar[r]\ar@{=}[d]&U_i[n]\ar[r]^-{\epsilon_{i,j}[n]}\ar@{_{(}->}[d]&
U_j[n-(j-i)]\ar[r]\ar@{_{(}->}[d]& 0\\
0\ar[r]& U_i[j-i]\ar[r]& U_i[n+1]\ar[r]^-{\epsilon_{i,j}[n+1]}&
U_j[n-(j-i)+1]\ar[r]& 0 ,}$$ where the map ${\epsilon_{i,j}[n]}$ is
induced by the canonical inclusion $U_i[j-i]\hookrightarrow U_i[n]$.
Thus, by taking the direct limit of the above diagram, we obtain the
following canonical exact sequence
$$(*)\quad 0\lra U_i[j-i]\lra U_i[\infty]\lraf{\epsilon_{i,j}}U_{j}[\infty]\lra 0, $$
where
$\epsilon_{i,j}:=\displaystyle\varinjlim_{n}\epsilon_{i,j}[n]$. This
finishes the proof of the first assertion in $(2)$. In the
following, we shall show that $\epsilon_{i,j}=\epsilon_{i+m,\,j+m}$
and $\epsilon_{i,j}\,\epsilon_{j,p}=\epsilon_{i,p}$ for any $p>j$.
In fact, the former clearly follows from
$\epsilon_{i,j}[n]=\epsilon_{i+m,\,j+m}[n]$  for any $n>j-i$, since
$U_i=U_{i+m}$  and  $U_j=U_{j+m}$  by our convention. As for the
latter, one can check  that, for any $u>p-i$, the composition of
$$\epsilon_{i,j}[u]: U_i[u]\lra U_j[u-(j-i)]\quad \mbox{and}\quad \epsilon_{j,
p}[u-(j-i)]: U_j[u-(j-i)]\lra U_p[u-(p-i)]$$ coincides with
$\epsilon_{i, p}[u]: U_i[u]\lra U_p[u-(p-i)]$. So, we have
$\epsilon_{i,j}[u]\,\epsilon_{j, p}[u-(j-i)]=\epsilon_{i, p}[u]$.
Consequently, by taking the direct limit of the two-sides of the
equality, we have $\epsilon_{i,j}\,\epsilon_{j,p}=\epsilon_{i,p}$
for any $p>j$. This completes the proof of $(2)$.

$(3)$ If  we apply $\Hom_R(U_i[\infty],-)$ to the sequence $(*)$ in
the proof of $(2)$, then we can get the following exact sequence:
$$0\ra \Hom_R(U_i[\infty],U_i[j-i])\ra \Hom_R(U_i[\infty],U_i[\infty])
\lraf{(\epsilon_{i,j})^*} \Hom_R(U_i[\infty],U_{j}[\infty])\ra
\Ext^1_R(U_i[\infty],U_i[j-i]).$$ Note that
$\Hom_R(U_i[\infty],U_i[j-i])=0$ by $(1)$. Thus, to prove that
$(\epsilon_{i,j})^*$ is an isomorphism, it suffices to show
$\Ext^1_R(U_i[\infty],U_i[j-i])=0$. In fact, this follows from
$\Ext^1_R(U_i[\infty],U_i[j-i])\simeq
D\Hom_R\big(\tau^-(U_{i}[j-i]),U_i[\infty]\big) \simeq
D\Hom_R(U_{i+1}[j-i],U_i[\infty]=0,$  where the last equality holds
for $1\le j-i<m$ by (1).

Next, if we apply $\Hom_R(-, U_j[\infty])$ to the sequence ($*$),
then we get the following exact sequence:
$$0\ra \End_R(U_j[\infty])\lraf{(\epsilon_{i,j})_*}
\Hom_R(U_i[\infty],U_i[\infty])\lra \Hom_R(U_i[j-i],U_j[\infty]).$$
Since $1\leq j-i <m$, we have $\Hom_R(U_i[j-i],U_j[\infty])=0$, and
therefore $(\epsilon_{i,j})_*$ is an isomorphism.

Now, it follows from the isomorphisms $(\epsilon_{i,j})^*$ and
$(\epsilon_{i,j})_*$ that  the map $$\varphi_{i,j}:
\End_R(U_i[\infty])\ra \End_R(U_j[\infty])$$ in $(3)$ is
well-defined and thus a ring isomorphism.

$(4)$ By definition, for $1\leq r,\,s,\, t\leq m$, one can check
$$\pi_{r,s}\,\pi_{s,t}=\epsilon_{r, s+\Delta_{r,s}\,m}\,\epsilon_{s, t+\Delta_{s,t}\,m}
=\epsilon_{r, s+\Delta_{r,s}\,m}\,\epsilon_{s+\Delta_{r,s}\,m,\,
t+(\Delta_{s,t}+\Delta_{r,s})\,m}=\epsilon_{r,\, t+(\Delta_{r,s} +
\Delta_{s,t})\,m}.
$$ On the one hand,  for
any $p>r$ and $q>r$, we infer from $(2)$ that
$\epsilon_{r,p}=\epsilon_{r,q}$ if and only if $p=q$. On the other
hand, we always have
$\Delta_{r,s}+\Delta_{s,t}-\Delta_{r,t}\in\{0,1\}$. Consequently,
the first statement in (4) follows. In particular, this implies that
$\pi_{i,j}\pi_{j,j}=\pi_{i,i}\pi_{i,j}$ for $1\leq i<j\leq m$. By
the definition of $\varphi_{i,j}$ in (3), we can prove the second
statement in (4).

$(5)$ Set $D_i:=\End_R(U_i[\infty])$. It follows from \cite[Section
4.4]{R} that $D_i$ is a complete discrete valuation ring. Let
$\mathfrak{m}$ be the unique maximal ideal of $D_i$. We shall prove
that $\pi_{i.i}$ is a prime element of $D_i$, that is, $\mathfrak
m=\pi_{i,i}D_i=D_i\pi_{i,i}$. Indeed, by applying $\Hom_R(-,
U_i[\infty])$ to the following exact sequence:
$$0\lra U_i[m]\lra
U_i[\infty]\lraf{\pi_{i,i}} U_{i}[\infty]\lra 0, $$ we obtain
another exact sequence of right $D_i$-modules:
$$0\lra D_i\lraf{(\pi_{i,i})_*}
D_i\lra\Hom_R(U_i[m], U_{i}[\infty])\lra 0, $$ due to
$\Ext_R^1(U_i[\infty], U_i[\infty])=0$, which follows from
\cite[Section 4.5]{R}. To show $\mathfrak m=\pi_{i,i}D_i$, we first
claim that $\Hom_R(U_i[m], U_{i}[\infty])\simeq \Hom_R(U_i,
U_{i}[\infty])\simeq D_i/{\mathfrak m}$ as right $D_i$-modules.

Let $$0\lra U_i\lra U_i[m]\lraf{\epsilon_{i,i+1}[m]}
U_{i+1}[m-1]\lra 0$$ be the exact sequence defined in (2). Then we
get the following exact sequence of $k$-modules:
$$
\Hom_R(U_{i+1}[m-1],U_i[\infty])\lra \Hom_R(U_i[m],U_i[\infty])\lra
\Hom_R(U_i,U_i[\infty])\lra \Ext_R^1(U_{i+1}[m-1],U_i[\infty]).$$
Since
$\Hom_R(U_{i+1}[m-1],U_i[\infty])=0=\Ext_R^1(U_{i+1}[m-1],U_i[\infty])$
by (1), we have $\Hom_R(U_i[m],U_i[\infty])\simeq
\Hom_R(U_i,U_i[\infty])$ as right $D_i$-modules.

It remains to show $\Hom_R(U_i, U_{i}[\infty])\simeq D_i/{\mathfrak
m}$ as right $D_i$-modules.  Let $$0\lra U_i\lraf{\zeta}
U_i[\infty]\lraf{\epsilon_{i,i+1}} U_{i+1}[\infty]\lra 0$$ be the
exact sequence defined in (2) with $\zeta$ the canonical inclusion.
Since $\Ext_R^1((U_{i+1})[\infty], U_i[\infty])=0$ by \cite[Section
4.5]{R}, we infer that, for any $f:U_i\to U_i[\infty]$, there is
$g\in D_i$ such that $f=\zeta g$. This means
$\Hom_R(U_i,U_i[\infty])=\zeta D_i$. Clearly, $\zeta D_i\simeq D_i/N
$ as right $D_i$-modules, where $N:=\{h\in D_i\mid \zeta h=0\}$. As
the canonical ring homomorphism from $D_i$ to $\End_R(U_i)$ via the
map $\zeta$ induces a ring isomorphism from $D_i/\mathfrak m$ to
$\End_R(U_i)$, we have $\zeta\mathfrak m=0$, that is, $\mathfrak
m\subseteq N$.  Since $D_i$ is a local ring and $N\subsetneq D_i$,
we get $N=\mathfrak m$, and therefore $\Hom_R(U_i,U_i[\infty])\simeq
D_i/\mathfrak m$ as right $D_i$-modules. This finishes the claim.

From the above claim, we conclude that $\mathfrak m$ coincides with
the image of $(\pi_{i,i})_*$, that is, $\mathfrak m=\pi_{i,i}D_i$.
Similarly, we can prove $\mathfrak m=D_i\pi_{i,i}$. This means that
$\pi_{i.i}$ is a prime element of $D_i$. As for the second statement
in $(5)$, we note that, for any $p\in\mathbb{N}$ and $1\leq q<m$,
the canonical inclusion map $U_i[pm+q]\ra U_i[pm+q+1]$ induces  an
isomorphism:
$$\Hom_R(U_i[pm+q+1],U_i[\infty])\lraf{\simeq}\Hom_R(U_i[pm+q],U_i[\infty]).$$
Consequently, we have the following isomorphisms of abelian groups:
$$\begin{array}{rl}
D_i&=\Hom_R\big(\displaystyle\varinjlim_n\,U_i[n],
U_i[\infty]\big)\simeq \displaystyle\varprojlim_n\,\Hom_R\big(U_i[n],U_i[\infty]\big) \\
&\simeq\displaystyle\varprojlim_n\,\Hom_R\big(U_i[(n-1)m+1],U_i[\infty]\big)\simeq
\displaystyle\varprojlim_n\, k[x]/(x^{n})\simeq k[[x]].\\
\end{array}$$
Here we need the assumption that $k$ is algebraically closed field.
Now, one can check directly that the composition of the above
isomorphisms yields a ring isomorphism $\varphi_i:D_i\to k[[x]]$,
which sends $\pi_{i,i}$ to $x$. This finishes the proof. $\square$

\bigskip
By Lemma \ref{lem3.3}(3), the rings $\End_R(U_i[\infty])$, with
$1\leq i\leq m$, are all isomorphic. From now on, we always identify
these rings, and simply denote them by $D({\mathcal C})$. Further,
we write ${\mathfrak m}({\mathcal C}) $ and $Q({\mathcal C})$ for
the maximal ideal of $D({\mathcal C})$ and the division ring of
fractions of $D({\mathcal C})$, respectively. In particular,
${\mathfrak m}({\mathcal C})=\pi_{i,i}D({\mathcal C})=D({\mathcal
C})\pi_{i,i}$.

\medskip
 Suppose
that $C$ is a $\mathbb{Z}$-module and $c\in C$. For $1\leq i, j\leq
m$,  we denote by $E_{i,j}(c)$ the $ m\times m$ matrix which has the
$(i,j)$-entry $c$, and the other entries $0$. For simplicity, we
write $E_{i,j}$ for $E_{i,j}(1)$ if $C$ is a ring with the identity
$1$.

\begin{Lem} For $1\leq i, j\leq m$, let $\pi_{i,j}$ be the homomorphisms defined
in  Lemma \ref{lem3.3}(4). Then there is a ring isomorphism
$$
\rho: \;\End_R(\bigoplus^m_{i=1}U_i[\infty])\lra\Gamma
\big({\mathcal C}\big):=\begin{pmatrix}
D({\mathcal C})             & D({\mathcal C})            & \cdots & D({\mathcal C}) \\
\mathfrak m({\mathcal C})   &  D({\mathcal C})            & \ddots &\vdots  \\
\vdots &  \ddots  &\ddots  & D({\mathcal C})\\
\mathfrak m({\mathcal C})  & \cdots    &  \mathfrak m({\mathcal C})   & D({\mathcal C}) \\
\end{pmatrix}_{m\times m}$$
which sends $E_{m,1}(\pi_{m,1})$ to $E_{m, 1}(\pi_{m,m})$ and
$E_{r,r+1}(\pi_{r,r+1})$ to $E_{r,r+1}$ for $1\leq r<m$, where the maximal ideal
${\mathfrak m}({\mathcal C})$ of the ring $D(\mathcal C)$ is generated by the element
$\pi_{m,m}$. \label{lem3.4}
\end{Lem}
{\it Proof.} For any $1\leq i<m$, by Lemma \ref{lem3.3}(2) and
$(4)$, we have the following exact sequence of $R$-modules:
$$0\lra U_i[m-i]\lra U_i[\infty]\lraf{\pi_{i,m}} U_{m}[\infty]\lra
0.$$ Summing up these sequences,  we  can get the following exact
sequence:
$$0\lra \bigoplus^{m-1}_{i=1}U_i[m-i]\lra \bigoplus^{m}_{j=1}U_j[\infty]
\lraf{\xi} U_{m}[\infty]^{(m)}\lra 0,$$ where
$\xi:=\mbox{diag}\big(\pi_{1,m},\,\pi_{2,m},\,\cdots,\,\pi_{m-1,m},\,1\big)$,
the $m\times m$ diagonal matrix with $\pi_{i,m}$ in the
$(i,i)$-position for $1\leq i<m$, and with 1 in the
$(m,m)$-position.

Let $D:=\End_R(U_m[\infty])$,  and let $\mathfrak m$ be the unique
maximal ideal of $D$. Set
$\Lambda:=\End_R(\displaystyle\bigoplus^m_{j=1}U_j[\infty])$. Since
$\Hom_R\big(U_i[m-i],\,U_m[\infty]\big)=0$ for $1\leq i<m$, we see
that, for any $g\in\Lambda$, there exists a unique homomorphism $f$
and a unique homomorphism $h$ such that the following diagram is
commutative:
$$\xymatrix{
 0\ar[r]& \displaystyle\bigoplus^{m-1}_{i=1}U_i[m-i] \ar[r]\ar@{-->}_-{f}[d]&
 \displaystyle\bigoplus^{m}_{j=1}U_j[\infty]\ar[r]^-{\xi}\ar[d]_-{g}&U_{m}[\infty]^{(m)}\ar[r]\ar@{-->}_-{h}[d]&0\\
0\ar[r]&  \displaystyle\bigoplus^{m-1}_{i=1}U_i[m-i]\ar[r]&
 \displaystyle\bigoplus^{m}_{j=1}U_j[\infty]\ar[r]^-{\xi}&U_{m}[\infty]^{(m)}\ar[r]&0.
}
$$
This yields a ring homomorphism $\rho:\Lambda\to M_m(D)$ defined by
$g\mapsto h$. More precisely, if $g=\big(g_{u,v}\big)_{1\leq u,
v\leq m}\in \Lambda$ with $g_{u,v}\in\Hom_R(U_u[\infty],
U_v[\infty])$, then $h=\big(h_{u,v}\big)_{1\leq u, v\leq m}\in
M_m(D)$ with $h_{u,v}\in D$ satisfying

(a) $g_{u,v}\pi_{v,m}=\pi_{u,m}h_{u,v}$ if $u<m$ and $v<m$,

(b) $h_{m,v}=g_{m,v}\pi_{v,m}$ if $u=m$ and $v<m$,

(c) $g_{u,m}=\pi_{u,m}h_{u,m}$ if $u<m$ and $v=m$,  and

(d) $h_{m,m}=g_{m,m}$.

\noindent In particular, the map $\rho$ sends $E_{u,u}$ in $\Lambda$
to $E_{u,u}$ in $M_m(D)$. In this sense, we may write
$\rho=\big(\rho_{u,v}\big)_{1\leq u, v\leq m}$, where
$\rho_{u,v}:\Hom_R(U_u[\infty], U_v[\infty])\to D$ is defined by
$g_{u,v}\mapsto h_{u,v}$.

Clearly, $\rho$ is injective since $\Hom_R(U_j[\infty],U_i[m-i])=0$
for  $1\leq j\leq m$ and $1\leq i<m$ by Lemma \ref{lem3.3}(1). In
the following, we shall determine the image of $\rho$, which is
clearly a subring of $M_m(D)$.

On the one hand, for any $a\in\End_R(U_u[\infty])$,
$b\in\Hom_R(U_u[\infty], U_v[\infty])$ and
$c\in\End_R(U_v[\infty])$, we have
$(abc)\rho_{u,v}=(a)\rho_{u,u}(b)\rho_{u,v}(c)\rho_{v,v}$. On the
other hand, it follows from Lemma \ref{lem3.3}(3) that $\rho_{u,u}$
is always  a ring isomorphism, and the left
$\End_R(U_u[\infty])$-module $\Hom_R(U_u[\infty],U_v[\infty])$ is
freely generated by $\pi_{u,v}$ for  $1\leq u\neq v\leq m$. This
implies that the image of $\rho$ coincides with the $m\times m$
matrix ring having $D\,(\pi_{u,v})\rho_{u,v}$ in the
$(u,v)$-position if $1\leq u\neq v\leq m$, and $D$ otherwise.
Moreover, by Lemmata \ref{lem3.3}(3) and (4), if $1\leq s<t<m$ and
$1\leq w<m$, we can form the following commutative diagrams:
$$
\xymatrix{U_s[\infty]\ar[r]^-{\pi_{s,m}}\ar[d]^-{\pi_{s, t}}&
U_m[\infty]\ar@{=}[d]\\
U_t[\infty]\ar[r]^-{\pi_{t,m}}& U_m[\infty] ,}\quad\,
\xymatrix{U_t[\infty]\ar[r]^-{\pi_{t,m}}\ar[d]^-{\pi_{t, s}}&
U_m[\infty]\ar[d]^-{\pi_{m,m}}\\
U_s[\infty]\ar[r]^-{\pi_{s,m}}& U_m[\infty]
,}\quad\,\xymatrix{U_m[\infty]\ar@{=}[r]\ar[d]^-{\pi_{m, w}}&
U_m[\infty]\ar[d]^-{\pi_{m,m}}\\
U_w[\infty]\ar[r]^-{\pi_{w,m}}& U_m[\infty] ,}\quad\,
\xymatrix{U_w[\infty]\ar[r]^-{\pi_{w,m}}\ar[d]^-{\pi_{w, m}}&
U_m[\infty]\ar@{=}[d]\\
U_m[\infty]\ar@{=}[r]& U_m[\infty].}
$$
In other words, we have $(\pi_{s,t})\rho_{s, t}=1=(\pi_{w,
m})\rho_{w, m}$ and $(\pi_{t,s})\rho_{t,s}=\pi_{m, m}=(\pi_{m,
w})\rho_{m, w}$. Thus, the image of $\rho$ is equal to the $m\times
m$ matrix ring having $D\,\pi_{m, m}$ as the $(p,q)$-entry for
$1\leq q<p\leq m$, and $D$ as the other entries. By Lemma
\ref{lem3.3}(5), we know $\mathfrak m=D\,\pi_{m.m}$. Now, by
identifying $D$ with $D({\mathcal C})$ and $\mathfrak m$ with
$\mathfrak m ({\mathcal C})$, we infer that the image of $\sigma$
coincides with the ring $\Gamma \big({\mathcal C}\big)$ defined in
Lemma \ref{lem3.4}. Therefore, we conclude that $\rho:\Lambda\to
\Gamma \big({\mathcal C}\big)$ is a ring isomorphism which sends
$E_{m,1}(\pi_{m,1})$ to $E_{m, 1}(\pi_{m,m})$ and $E_{r,
r+1}(\pi_{r,r+1})$ to $E_{r, r+1}$ for $1\leq r<m$. This completes
the proof. $\square$

\medskip
Combining Lemma \ref{lem3.3}(5) with Lemma \ref{lem3.4}, we then
obtain the following result which will be used for the calculation
of stratifications of derived module categories in the next section.

\begin{Koro}
For $1\leq i, j\leq m$, let $\pi_{i,j}$ be the homomorphisms defined
in  Lemma \ref{lem3.3}(4). Assume that $k$ is an algebraically
closed field.  Then there exists a ring isomorphism
$$
\sigma: \;\End_R(\bigoplus^m_{i=1}U_i[\infty])\lra\Gamma
(m):=\begin{pmatrix}
k[[x]] &  k[[x]]  & \cdots & k[[x]]\\
(x)    &  k[[x]]  & \ddots & \vdots\\
\vdots &  \ddots  &\ddots  & k[[x]]\\
(x)    &  \cdots        &(x)     & k[[x]]\\
\end{pmatrix}_{m\times m}$$
which sends $E_{m,1}(\pi_{m,1})$ to $E_{m, 1}(x)$ and
$E_{r,r+1}(\pi_{r,r+1})$ to $E_{r,r+1}$ for $1\leq r<m$.
\label{cor3.4}
\end{Koro}

\subsection{Universal localizations at simple regular modules }

 From now on, let us fix a non-empty subset
$\mathcal{U}$ of ${\mathscr S}$, where ${\mathscr S}$ is a complete
set of isomorphism classes of all simple regular $R$-modules. Denote
by $\lambda:R\to R_\mathcal{U}$ the universal localization of $R$ at
$\mathcal{U}$. It follows from \cite[Theorems 4.9, 5.1, and
5.3]{Sch1} that $\lambda$ is injective and $R_{\mathcal{U}}$ is
hereditary. Moreover, it is shown in \cite[Corollary 4.6(2),
4.7]{HJ} and \cite{HJ2} that the $R$-module
$$
T_\mathcal{U}:=R_{\mathcal{U}}\oplus R_{\mathcal{U}}/R
$$
is a tilting module with $\Hom_R(R_\mathcal{U}/R, R_\mathcal{U})=0$.

Suppose
$$(*)\quad 0\lra R\lraf{\lambda} R_\mathcal{U}\lraf{\pi} R_\mathcal{U}/R\lra 0,$$
is the canonical exact sequence of $R$-modules with $\pi$ the
canonical surjection. Set $B:=\End_R(T_\mathcal{U})$,
$S:=\End_R(R_\mathcal{U}/R)$ and $\Sigma:=\{S\otimes_R f_U\mid
U\in\mathcal{U}\}$. Recall that the right multiplication map
$\mu:R\to S$ defined by $r\mapsto(y\mapsto yr)$ for $r\in R$ and
$y\in S/R$, is a ring homomorphism, which endows $S$ with a natural
$R$-$R$-bimodule structure.

Let $\mathcal{U}^+$ be the full subcategory of $R$-Mod, defined by
$$\mathcal{U}^+:=\{X\in R\Modcat\mid\Ext^i_R(U,X)=0 \mbox{ for\, all
\,} U\in\mathcal{U}\; \mbox{and\,all \,}{i\in\mathbb{N}}\}.$$ For
example, the Pr\"ufer module $V[\infty]$ for $V\in {\mathscr
S}\setminus \mathcal{U}$ lies in ${\mathcal U}^+$ by Lemma
\ref{lem3.3}(1).

This subcategory has the following characterization, due to
\cite[Proposition 3.8]{HJ}.

\begin{Lem}
$\mathcal{U}^+$ coincides with the image of the restriction functor
$\lambda_*:R_\mathcal{U}\Modcat\to R\Modcat$. In particular, for any
$Y\in\mathcal{U}^+$,  the unit adjunction $\eta_Y: Y\to
R_\mathcal{U}\otimes_RY$, defined by $y\mapsto 1\otimes y$ for $y\in
Y$, is an isomorphism of $R$-modules. \label{orth}
\end{Lem}

Thus, for an $R$-module $Y\in \mathcal{U}^+$, we may endow it with
an $R_{\mathcal U}$-module structure via the isomorphism $\eta_Y$,
and in this way, we consider the $R$-module $Y$ as an $R_{\mathcal
U}$-module. Note that this $R_{\mathcal U}$-module structure on $Y$
extended from the $R$-module structure is unique.

Concerning the universal localization $R_{\mathcal U}$ of $R$ at
$\mathcal U$, we have the following facts (see \cite[Proposition
1.10]{HJ2}, \cite{ Sch1} and \cite{CB}).

\begin{Lem}
$(1)$ Suppose that $\mathcal{U}$ contains no cliques. Then
$R_\mathcal{U}$ is a finite-dimensional  tame hereditary
$k$-algebra. In particular, the tilting $R$-module $T_\mathcal{U}$
is classical. Moreover, $\{R_\mathcal{U}\otimes_R V\mid V\in
{\mathscr S}\backslash \mathcal{U}\}$ is a complete set of
non-isomorphic simple regular $R_\mathcal{U}$-modules, and
$(R_\mathcal{U}\otimes_R V)[\infty]\simeq V[\infty]$ as $R_{\mathcal
U}$-modules for each $V\in{\mathscr S}\backslash \mathcal{U}$.

$(2)$ Suppose that $\mathcal{U}$ contains cliques. Then
$R_\mathcal{U}$ is a hereditary order. Moreover,
$\{R_\mathcal{U}\otimes_R V\mid V\in {\mathscr S}\backslash
\mathcal{U}\}$ is a complete set of non-isomorphic simple
$R_\mathcal{U}$-modules, and the injective envelope of the
$R_{\mathcal U}$-module $R_\mathcal{U}\otimes_R V$ is isomorphic to
$V[\infty]$ for each $ V\in {\mathscr S}\backslash \mathcal{U}$.

$(3)$ Suppose  $\mathcal{V}\subseteq {\mathscr S}\backslash
\mathcal{U}$. Then  $R_{\mathcal{U}\cup
\mathcal{V}}=(R_\mathcal{U})_{\ol{\mathcal{V}}}$, where
$\ol{\mathcal{V}}:=\{R_\mathcal{U}\otimes_R V\mid
V\in\mathcal{V}\}$. In particular, there are injective ring
epimorphisms $R_\mathcal{U}\lra R_{\mathcal{U}\cup \mathcal{V}}$ and
$R_{\mathcal{U}\cup \mathcal{V}}\lra R_{\mathscr S}$. \label{lem3.1}
\end{Lem}

As remarked in \cite{CB}, in the case of Lemma \ref{lem3.1}(1), the
set of simple regular $R_{\mathcal{U}}$-modules in a clique is of
the form
$$ \{R_{\mathcal U}\otimes_RV\mid V\in {\mathcal C}, V\not\in
\mathcal{U}\},$$ where $\mathcal C$ is a clique of $R$. Further, by
Lemma \ref{lem3.1}(1), for each $V\in{\mathcal C}\backslash
\mathcal{U}$, the Pr\"ufer modules corresponding to $R_{\mathcal
U}\otimes_RV$ and  to $V$ are isomorphic.  In particular, they have
the isomorphic endomorphism ring.

Thus, if $\mathcal{C}_1, \mathcal{C}_2, \cdots, \mbox{and}\,
\mathcal{C}_s$ are all cliques from non-homogeneous tubes and if
$\mathcal{U}$ is a union of $c(\mathcal{C}_i)-1$ simple regular
$R$-modules from each $\mathcal{C}_i$, then each clique of
$R_{\mathcal U}$ consists of only one single element. This implies
that $R_\mathcal{U}$ has only two isomorphism classes of simple
modules. If, in addition, the field $k$ is algebraically closed,
then $R_\mathcal{U}$ is Morita equivalent to the Kronecker algebra.
In this case, since the set of cliques of the Kronecker algebra are
parameterized by $\mathbb{P}^1(k)$, we see that the set of cliques
of an arbitrary tame hereditary $k$-algebra can be indexed by
$\mathbb{P}^1(k)$.

A description  of the structure of the module $R_{\mathcal U}/R$ was
first given in \cite{Sch2}, and a further substantial discussion is
carried out recently in \cite{HJ2}. Especially, the following lemma
is proved in \cite{HJ2}, where the field $k$ is required to be
algebraically closed. In fact, one can check that, if $k$ is an
arbitrary field, all of the arguments in the proof of the lemma in
\cite{HJ2} are still valid except some mild changes. For instance,
the field $k$ should be replaced by certain division rings in most
of the proofs.

\begin{Lem}\label{moddecom}
$(1)$ The $R$-module $R_\mathcal{U}/R$ is a direct union of finite
extensions of modules in $\mathcal{U}$.

$(2)$ Let ${\bf t}\subset R\modcat$ be a tube of rank $m>1$, and let
$\mathcal{U}=\{U_1,U_2,\cdots, U_{m-1}\}$ be a set of $m-1$ simple
regular modules in $\bf{t}$ such that $U_{i+1}=\tau^-U_i$ for all
$1\leq i\le m-1.$ Then
$$R_\mathcal{U}/R\simeq U_1[m-1]^{(\delta_{U_1})}\oplus U_2[m-2]^{(\delta_{U_2})}
\oplus \cdots \oplus U_{m-1}[1]^{(\delta_{U_{m-1}})},$$ with
$\delta_{U_j}:=\dim_{\End_R(U_j)}\Ext^1_R(U_j,R)$ for $1\leq j\leq
m-1$. Moreover, $R_\mathcal{U}\otimes_R U_m\simeq U_m[m]$ as
$R_{\mathcal U}$-modules.

$(3)$ If $\mathcal{U}$ is a union of cliques, then, for any finitely
generated projective $R$-module $P$,
$$_R(R_\mathcal{U}/R)\otimes_R P\simeq \bigoplus_{U\in\mathcal{U}}U[\infty]^{(\delta_{U,P})},$$
where $\delta_{U,P}:=\dim_{\End_R(U)}\Ext^1_R(U,P).$
\end{Lem}
Next, we shall show that $R_{\mathcal U}$ and
 $\End_R(R_{\mathcal U}/R)$ can be interpreted as the tensor product
and direct sum of some rings, respectively.

\begin{Lem}\label{3.4}
Let $\mathcal{U}=\mathcal{U}_0\dot{\cup}
\mathcal{U}_1\subseteq{\mathscr S}$ such that ${\mathcal U}_0$
contains no cliques and ${\mathcal U}_1$ is a union of cliques. Then
the following statements are true:

$(1)$ $\mathcal{U}_0\subseteq \mathcal{U}_1^{+}$,
$\mathcal{U}_1\subseteq \mathcal{U}_0^{+}$, $R_{\mathcal{U}}\simeq
R_{\mathcal{U}_1}\otimes_R R_{\mathcal{U}_0}$ as
$R_{\mathcal{U}_1}$-$R_{\mathcal{U}_0}$-bimodules, and
$R_{\mathcal{U}}/R_{\mathcal{U}_1}\simeq R_{\mathcal{U}_1}\otimes_R
(R_{\mathcal{U}_0}/R)$ as $R_{\mathcal{U}_1}$-$R$-bimodules.

$(2)$ There is a ring isomorphism
$$\varphi:\End_R(R_{\mathcal{U}}/R)\lra\End_R(R_{\mathcal{U}_0}/R)\times
\End_{R_{\mathcal{U}_0}}(R_\mathcal{U}/R_{\mathcal{U}_0}).$$

\end{Lem}
{\it Proof}. $(1)$ By the assumption on $\mathcal U$, if
$U\in{\mathcal U}_0$ and $V\in {\mathcal U}_1$ then they belong to
different tubes, and therefore $\mathcal{U}_0\subseteq
\mathcal{U}_1^{+}$ and $\mathcal{U}_1\subseteq \mathcal{U}_0^{+}$.

By Lemma \ref{orth}, the unit adjunction $\eta_{U}: U\to
R_{\mathcal{U}_0}\otimes_R U$ is an isomorphism of $R$-modules for
any $U\in\mathcal{U}_1$. This implies that every module in
$\mathcal{U}_1$ can be endowed with a unique
$R_{\mathcal{U}_0}$-module structure that preserves the given
$R$-module structure via the universal localization
$\lambda_{0}:R\to R_{\mathcal{U}_0}$. Consequently, it follows from
Lemma \ref{lem3.1}(3) that
$R_\mathcal{U}=(R_{\mathcal{U}_0})_{\mathcal{U}_1}$. Moreover, we
can construct the following exact commutative diagram of
$R$-modules:
\begin{eqnarray*} \xymatrix{
& &       & 0\ar[d]             &0\ar[d]         &\\
&
0\ar[r]&R\ar[r]^-{\lambda_0}\ar@{=}[d]&R_{\mathcal{U}_0}\ar[d]^-{\lambda_1}
\ar[r]^-{\pi_0}\ar[d]&R_{\mathcal{U}_0}/R\ar[r]\ar[d]^-{\lambda_2}&0\\
(*) & 0\ar[r]& R\ar[r]^-{\lambda}
&R_\mathcal{U}\ar[d]^-{\pi_1}\ar[r]^-{\pi}&
R_\mathcal{U}/R\ar[r]\ar[d]^-{\pi_2}\ar[d]&0\\
& &        & R_\mathcal{U}/R_{\mathcal{U}_0}\ar@{=}[r]\ar[d] &
R_\mathcal{U}/R_{\mathcal{U}_0}\ar[d]& \\
& &        &                     0 &  0,& \\
}
\end{eqnarray*}
where  $\lambda_1$ is the universal localization of
$R_{\mathcal{U}_0}$ at $\mathcal{U}_1$, and $\lambda_2$ is the
canonical injection induced by $\lambda_1$, and where $\pi_0$,
$\pi_1$ and $\pi_2$ are canonical surjections.

Clearly, $R_{\mathcal{U}_0}$ is a finite-dimensional tame hereditary
algebra by Lemma \ref{lem3.1}(1). From
$R_\mathcal{U}=(R_{\mathcal{U}_0})_{\mathcal{U}_1}$ we see that
$R_{\mathcal{U}}/R_{\mathcal{U}_0}$ is a direct union of finite
extensions of modules in $\mathcal{U}_1$ by Lemma \ref{moddecom}(1).
Since $R_{\mathcal{U}_1}$ is the universal localization of $R$ at
$\mathcal{U}_1$, we have $\Tor^R_i(R_{\mathcal{U}_1}, V)=0$ for any
$i\geq 0$ and $V\in\mathcal{U}_1$. Note that the $i$-th left derived
functor $\Tor^R_i(R_{\mathcal{U}_1},-):R\Modcat\to
\mathbb{Z}\Modcat$ commutes with direct limits.  Thus
$\Tor^R_i(R_{\mathcal{U}_1}, R_{\mathcal{U}}/R_{\mathcal{U}_0})=0$
for any $i\ge 0$, which implies that the homomorphisms
$R_{\mathcal{U}_1}\otimes_R\lambda_1$ and
$R_{\mathcal{U}_1}\otimes_R\lambda_2$ are isomorphisms. Moreover, by
Lemma \ref{lem2.5}, we have
$R_\mathcal{U}=(R_{\mathcal{U}_1})_{\ol{\mathcal{U}}_0}$ with
$\ol{\mathcal{U}}_0:=\{R_{\mathcal{U}_1}\otimes_R U\mid
U\in\mathcal{U}_0\}$, and therefore $R_\mathcal{U}$ can be regarded
as an $R_{\mathcal{U}_1}$-module. Consequently, the canonical
multiplication map $\nu_2:R_{\mathcal{U}_1}\otimes_R
R_\mathcal{U}\to R_\mathcal{U}$ is an isomorphism.

Now we apply the tensor functor $R_{\mathcal{U}_1}\otimes_R-$ to the
diagram $(*)$, and get the following exact commutative diagram of
$R_{\mathcal{U}_1}$-$R$-bimodules:
$$
\xymatrix{
&R_{\mathcal{U}_1}\otimes_RR\ar[rr]^-{R_{\mathcal{U}_1}\otimes_R\lambda_0}
\ar@{=}[d]&&R_{\mathcal{U}_1}\otimes_R
R_{\mathcal{U}_0}\ar[d]^-{R_{\mathcal{U}_1}\otimes_R\lambda_1}_-{\simeq}
\ar[rr]^-{R_{\mathcal{U}_1}\otimes_R\pi_0}\ar[d] &&
R_{\mathcal{U}_1}\otimes_R(R_{\mathcal{U}_0}/R)
\ar[r]\ar[d]^-{R_{\mathcal{U}_1}\otimes_R\lambda_2}_-{\simeq}&0\\
&R_{\mathcal{U}_1}\otimes_RR
\ar[d]_-{\simeq}^-{\nu_1}\ar[rr]^-{R_{\mathcal{U}_1}\otimes_R\lambda}
&&R_{\mathcal{U}_1}\otimes_R
R_\mathcal{U}\ar[d]_-{\simeq}^-{\nu_2}\ar[rr]^-{R_{\mathcal{U}_1}\otimes_R\pi}&&
R_{\mathcal{U}_1}\otimes_R (R_\mathcal{U}/R)\ar[r]\ar@{-->}[d]&0\\
0\ar[r]&R_{\mathcal{U}_1}\ar[rr]&& R_\mathcal{U}\ar[rr]&&
R_{\mathcal{U}}/R_{\mathcal{U}_1}\ar[r]&0,}
$$
where $\nu_1$ is the multiplication map. Thus $R_{\mathcal{U}}\simeq
R_{\mathcal{U}_1}\otimes_R R_{\mathcal{U}_0}$ as
$R_{\mathcal{U}_1}$-$R_{\mathcal{U}_0}$-bimodules, and
$R_{\mathcal{U}}/R_{\mathcal{U}_1}\simeq R_{\mathcal{U}_1}\otimes_R
(R_{\mathcal{U}_0}/R)$ as $R_{\mathcal{U}_1}$-$R$-bimodules.

$(2)$ Note that $\mathcal{U}_1\subseteq \mathcal{U}_0^{+}$ and
$\mathcal{U}_0\subseteq \mathcal{U}_1^{+}$, and that $\mathcal{U}_0$
and $\mathcal{U}_1$ consist of finitely presented modules of
projective dimension one.  By  Lemmata \ref{lem3.1}(1) and
\ref{moddecom}(1), we can write $R_{\mathcal{U}}/R_{\mathcal{U}_0}
=\displaystyle\varinjlim_{\alpha}X_{\alpha}$ with $X_{\alpha}\in
{\mathcal F}({\mathcal U}_1)$. Then, by Lemma \ref{2.7}, we have the
following isomorphisms:
$$(**)\quad \Ext^j_R(R_{\mathcal{U}_0}/R,
R_\mathcal{U}/R_{\mathcal{U}_0})\simeq\displaystyle\varinjlim_{\alpha}
\Ext_R^j(R_{\mathcal{U}_0}/R,
X_{\alpha})=0=\displaystyle\varprojlim_{\alpha}
\Ext_R^j(X_{\alpha},R_{\mathcal{U}_0}/R)\simeq\Ext_R^j(R_\mathcal{U}/R_{\mathcal{U}_0},
R_{\mathcal{U}_0}/R)$$ for any $j\geq 0$. Particularly, the
canonical exact sequence
$$0\lra R_{\mathcal{U}_0}/R\lraf{\lambda_2}
R_\mathcal{U}/R\lraf{\pi_2} R_\mathcal{U}/R_{\mathcal{U}_0}\lra 0$$
splits in $R\Modcat$, that is, $R_\mathcal{U}/R\simeq
R_{\mathcal{U}_0}/R \;\oplus\; R_\mathcal{U}/R_{\mathcal{U}_0}$ as
$R$-modules. Since $R\ra R_{{\mathcal U}_0}$ is a ring epimorphism,
we have $\End_R(R_{\mathcal U}/R_{{\mathcal
U}_0})=\End_{R_{{\mathcal U}_0}}(R_{\mathcal U}/R_{{{\mathcal
U}_0}})$. Thus it follows from ($**$) for $j=0$ that
$$\End_R(R_{\mathcal{U}}/R)\simeq \End_R(R_{\mathcal{U}_0}/R)\times
\End_{R_{\mathcal{U}_0}}(R_\mathcal{U}/R_{\mathcal{U}_0}).$$ This
completes the proof of (2). $\square$

\subsection{Proof of Theorem 1.1}

Before we start with the proof of the main result, Theorem
\ref{th1.1}, we have to make the following preparations.

\begin{Lem}
Let $\mathcal{U}=\mathcal{U}_0\dot{\cup}
\mathcal{U}_1\subseteq{\mathscr S}$ such that $\mathcal{U}_1$ is a
union of cliques and $\mathcal{U}_0$ does not contain any cliques.
Set
$\Lambda:=\End_{R_{\mathcal{U}_0}}(R_\mathcal{U}/R_{\mathcal{U}_0})$
and
$\Theta:=\{\Lambda\otimes_{R_{\mathcal{U}_0}}({R_{\mathcal{U}_0}}\otimes_Rf_V)\mid
V\in \mathcal{U}_1\}$, $S:=\End_R(R_{\mathcal U}/R)$ and $\Sigma
:=\{S\otimes_Rf_U\mid U\in {\mathcal U}\}$. Then $S_\Sigma$ is
isomorphic to the universal localization $\Lambda_\Theta$ of
$\Lambda$ at $\Theta$. \label{lem3.2}
\end{Lem}

{\it Proof.} By Lemma \ref{orth}, we have
$R_{\mathcal{U}_0}\otimes_R V\simeq V$ as $R$-modules for each
$V\in\mathcal{U}_1$. Combining this with Lemma \ref{lem3.1}(1), we
see that $\mathcal{U}_1$ can be seen as a set of simple regular
$R_{\mathcal{U}_0}$-modules, and therefore
$R_\mathcal{U}=(R_{\mathcal{U}_0})_{\mathcal{U}_1}$ by Lemma
\ref{lem3.1}(3).
 More precisely, for any $V\in
\mathcal{U}_1$, we fix a minimal projective presentation
$$
0\lra P_{1}\lraf{f_V} P_0\lra V\lra 0
$$ of $V$ in $R\modcat$,
and get a projective presentation of $V$ in
$R_{\mathcal{U}_0}\modcat:$
$$
0\lra R_{\mathcal{U}_0}\otimes_R
P_{1}\lraf{R_{\mathcal{U}_0}\otimes_Rf_V}
R_{\mathcal{U}_0}\otimes_RP_0\lra V\lra 0.$$ This is due to the fact
that $\Tor^R_1(R_{\mathcal{U}_0},V)\simeq\Tor^R_1(R_{\mathcal{U}_0},
R_{\mathcal{U}_0}\otimes_RV)\simeq\Tor^{R_{\mathcal{U}_0}}_1(R_{\mathcal{U}_0},
R_{\mathcal{U}_0}\otimes_RV)=0$. Thus, $R_\mathcal{U}$ is the
universal localization of $R_{\mathcal{U}_0}$ at the set
$\{R_{\mathcal{U}_0}\otimes_Rf_V\mid V\in\mathcal{U}_1\}$. Note that
$R_{\mathcal{U}_0}$ is a tame hereditary $k$-algebra  by Lemma
\ref{lem3.1}(1).

Let
$\Lambda:=\End_{R_{\mathcal{U}_0}}(R_\mathcal{U}/R_{\mathcal{U}_0})$
and
$\Theta:=\{\Lambda\otimes_{R_{\mathcal{U}_0}}({R_{\mathcal{U}_0}}\otimes_R
f_V)\mid V\in \mathcal{U}_1\}.$ In the following, we shall show that
$S_\Sigma$ is isomorphic to $\Lambda_\Theta$.

Let $\Gamma:=\End_R(R_{\mathcal{U}_0}/R)$ and
$\varphi=(\varphi_0,\varphi_1): S\ra\Gamma\times \Lambda$, where
$\varphi_0: S\to\Gamma$ and $\varphi_1: S\to\Lambda$ are the ring
homomorphisms given in Lemma \ref{3.4}(2). Recall that $\mu:R\to S$
is the right multiplication map.  Set
$\mu_0=\mu\,\varphi_0:R\to\Gamma$ and
$\mu_1=\mu\,\varphi_1:R\to\Lambda$. Clearly, both $\mu_0$ and
$\mu_1$ are ring homomorphisms, through which $\Lambda$ and $\Gamma$
have a right $R$-module structure, respectively. Now, we write
$\Sigma:=\{S\otimes_R f_U\mid U\in\mathcal{U}\}$ as $\Phi \times
\Psi$ with $\Phi:=\{\Gamma\otimes_R f_U\mid U\in\mathcal{U}\}$ and
$\Psi:=\{\Lambda\otimes_R f_U\mid U\in\mathcal{U}\}$. Consequently,
the ring isomorphism $\varphi$ implies that
$S_\Sigma\simeq\Gamma_{\Phi}\times \Lambda_{\Psi}$. To finish the
proof, it suffices to prove that $\Gamma_{\Phi}=0$ and
$\Lambda_{\Psi}\simeq\Lambda_\Theta$.

Indeed, we write  $\Phi=\Phi_0\cup\Phi_1$ with
$\Phi_0:=\{\Gamma\otimes_R f_U\mid U\in\mathcal{U}_0\}$ and
$\Phi_1:=\{\Gamma\otimes_R f_U\mid U\in\mathcal{U}_1\}$. Then, by
Lemma \ref{lem2.5}, we have $\Gamma_{\Phi}\simeq
(\Gamma_{\Phi_0})_{\ol{\Phi}_1}$, where
$\ol{\Phi}_1:=\{\Gamma_{\Phi_0}\otimes_Rf_U\mid
U\in\mathcal{U}_1\}$. To prove $\Gamma_{\Phi}=0$, it suffices to
prove $\Gamma_{\Phi_0}=0$. Consider the canonical exact sequence of
$R$-modules:
$$0\lra R\lraf{\lambda_0} R_{\mathcal{U}_0}\lraf{\pi_0}
R_{\mathcal{U}_0}/R\lra 0.$$ By Lemma \ref{lem3.1}(1), the module
$T_{\mathcal{U}_0}:=R_{\mathcal{U}_0}\oplus R_{\mathcal{U}_0}/R$ is
a classical tilting $R$-module, and therefore $\D{R}$ is triangle
equivalent to $\D{\End_R(T_{\mathcal{U}_0})}$ in the recollement of
$\D{R}$, $\D{\End_R(T_{\mathcal{U}_0})}$ and $\D{\Gamma_{\Phi_0}}$
by Proposition \ref{lem2.4}. Thus $\Gamma_{\Phi_0}=0$ and
$\Gamma_{\Phi}=0.$

It remains to show $\Lambda_{\Psi}\simeq\Lambda_\Theta$.  Let
$\mu_2: R_{\mathcal{U}_0}\to\Lambda$ be the right multiplication map
defined by $r\mapsto(x\mapsto xr)$ for $r\in R_{\mathcal{U}_0}$ and
$x\in R_\mathcal{U}/R_{\mathcal{U}_0}$. Then, along the diagram
$(*)$ in the proof of Lemma \ref{3.4}, one can check that the
following diagram of ring homomorphisms
$$\xymatrix{
R\ar[r]^-{\lambda_0}\ar[d]^-{\mu} &
R_{\mathcal{U}_0}\ar[d]^-{\mu_2}\\
S\ar[r]^-{\varphi_1}& \Lambda }
$$
commutes. Now, we write $\Psi=\Psi_0\cup\Psi_1$ with
$$\Psi_0:=\{\Lambda\otimes_R f_U\mid
U\in\mathcal{U}_0\}\quad\mbox{and}\quad \Psi_1:=\{\Lambda\otimes_R
f_V\mid V\in\mathcal{U}_1\},$$ and claim $\Lambda_{\Psi_0}=\Lambda$.
It suffices to show that $\Lambda\otimes_Rf_U$ is an isomorphism for
any $U\in\mathcal{U}_0$. However, this follows from
$\Lambda\otimes_Rf_U\simeq\Lambda\otimes_{R_{\mathcal{U}_0}}
(R_{\mathcal{U}_0}\otimes_R f_U)$ and  $R_{\mathcal{U}_0}\otimes_R
f_U$ being an isomorphism by the definition of universal
localizations. Hence $\Lambda_{\Psi_0}=\Lambda$.

Now, we have $\ol{\Psi}_1:=\{\Lambda_{\Psi_0}\otimes_{\Lambda}h\mid
h\in\Psi_1\}=\Psi_1$. It follows from  Lemma \ref{lem2.5} that
$\Lambda_{\Psi}\simeq (\Lambda_{\Psi_0})_{\ol{\Psi}_1}\simeq
\Lambda_{\Psi_1}$. Further, we have
$\Lambda\otimes_Rf_V\simeq\Lambda\otimes_{R_{\mathcal{U}_0}}
(R_{\mathcal{U}_0}\otimes_Rf_V)$ for any $V\in\mathcal{U}_1$. By
comparing the elements in $\Theta$ with the ones in $\Psi_1$, one
knows immediately that $\Lambda_{\Psi}\simeq\Lambda_\Theta$, and
therefore $S_\Sigma\simeq\Lambda_\Theta$, finishing the proof.
$\square$

\bigskip
Next, we shall show that the universal localizations of interest for
us take actually the form of ad\`ele rings in the algebraic number
theory \cite{ne}. Before stating the following lemma, we first
recall some notations.

\smallskip
Let $\mathcal{C}$ be  a clique  of $R\modcat$. Recall that
$D({\mathcal C})$ stands for the endomorphism ring of a Pr\"ufer
module $V[\infty]$ with $V\in\mathcal{C}$. Note that $D({\mathcal
C})$ is a discrete valuation ring with the division ring
$Q({\mathcal C})$ of fractions of $D({\mathcal C})$.  Clearly,
$Q({\mathcal C})$ contains $D({\mathcal C})$ as a subring.

\begin{Lem}\label{lem3.5}
Suppose that  $\mathcal{U}\subseteq{\mathscr S}$ is a union of
cliques, say ${\mathcal U}=\cup_{i\in I}{\mathcal C}_i $ with $I$ an
index set. Let $S:=\End_R(R_\mathcal{U}/R)$ and
$\Sigma:=\{S\otimes_R f_U\mid U\in\mathcal{U}\}$. Then $S_\Sigma$ is
Morita equivalent to the ad\`ele ring
$$\mathbb{A}_{\,\mathcal{U}}:=\bigg\{\big(f_{i}\big)_{i\in
I}\in\prod_ {i\in I}Q({\mathcal C}_i)\;\big|\;f_{i}\in D({\mathcal
C}_i)\; \mbox{\,for almost all }\, i\in I \bigg\}.$$
\end{Lem}

{\it Proof}. For any finitely generated projective $R$-module $P$,
we always have $S\otimes_RP\simeq\Hom_R(R_\mathcal{U}/R,
(R_\mathcal{U}/R)\otimes_RP)$ as $S$-modules. Thus,  we can rewrite
$\Sigma = \{\Hom_R\big(R_\mathcal{U}/R,\,
(R_\mathcal{U}/R)\otimes_Rf_V\big)\mid V\in\mathcal{U}\}$. The whole
proof of Lemma \ref{lem3.5} will be proceeded in three steps.

\textbf{Step $(1)$}.  We provide an alternative form of the
homomorphism $(R_\mathcal{U}/R)\otimes_Rf_V$ for any
$V\in\mathcal{U}$.

In fact, this procedure can be done for each clique $\mathcal C$ in
$\mathcal U$. Let us give the details: Fix a clique $\mathcal
{C}\subseteq \mathcal{U}$ and an element $U\in \mathcal C$, and
choose a projective resolution $ \quad 0\lra P_{1}\lraf{f_U} P_0\lra
U\lra 0 $ of $U$ in $R\modcat$, where $P_1$ and $P_0$ are finitely
generated projective $R$-modules. As $\lambda: R\to R_\mathcal{U}$
is the universal localization of $R$ at $\mathcal{U}$, we  know that
$R_\mathcal{U}\otimes_Rf_U:R_\mathcal{U}\otimes_RP_1\to
R_\mathcal{U}\otimes_RP_0$ is an isomorphism. This yields  the
following exact and commutative diagram of $R$-modules:
{\begin{eqnarray*} \xymatrix{
&           &                               &0\ar[d]            &\\
&  0\ar[d]  &                               &U\ar[d]^-{\psi}    &\\
0\ar[r]&P_1\ar[r]^-{\lambda\otimes_RP_1}\ar[d]_-{f_U}&R_\mathcal{U}\otimes_RP_1\ar[r]^-{\pi\otimes_R
P_1} \ar[d]^-{R_\mathcal{U}\otimes_Rf_U}_-{\simeq} &
(R_\mathcal{U}/R)\otimes_RP_1\ar[r]
\ar[d]^-{(R_\mathcal{U}/R)\otimes{_R}f_U}&0\\
0\ar[r]&P_0\ar[r]^-{\lambda\otimes_RP_0}\ar[d] &
R_\mathcal{U}\otimes_RP_0\ar[r]^-{\pi\otimes_R P_0}
& (R_\mathcal{U}/R)\otimes_RP_0\ar[r]\ar[d]   & 0\\
& U \ar[d] &                             & 0                  & \\
& 0.       &                             &                    &}
\end{eqnarray*}}
Consider the following  short exact sequence of $R$-modules:
$$\xymatrix{
(a)\quad 0\ar[r]& U\ar[r]^-{\psi}&
(R_\mathcal{U}/R)\otimes_RP_1\ar[rr]^-{(R_\mathcal{U}/R)\otimes{_R}f_U}
&&(R_\mathcal{U}/R)\otimes_RP_0\ar[r]&  \; _{}0.}$$  On the one
hand, by Lemma \ref{moddecom}(3), we have
$$(R_\mathcal{U}/R)\otimes_RP_1\simeq
\bigoplus_{i\in I}\;\bigoplus_{V\in {\mathcal
C}_i}V[\infty]^{(n_V)}$$ for some $n_V\in\mathbb{N}$, where $n_U$ is
non-zero since $U$ can be embedded into
$(R_\mathcal{U}/R)\otimes_RP_1$.
 On the other hand, for $W\in\mathcal{U}$, we have $\Hom_R(U,
W[\infty])$ = $0$ if $W\ncong U$, and  $\Hom_R(U, U[\infty])\simeq
\End_R(U)$. Now, let
$$0\lra U\lraf{\zeta_U} U[\infty]\lraf{\pi_U}
(\tau^-U)[\infty]\lra 0$$ be the canonical exact sequence defined in
Lemma \ref{lem3.3}(2), where $\zeta_U$ is the canonical inclusion.
Set $D:=\End_R(U[\infty])$. Then $D$ is a discrete valuation ring.
Particularly, it is a local ring with a unique maximal ideal
$\mathfrak m$. By the proof of Lemma \ref{lem3.3}(5), we know that
$\Hom_R(U, U[\infty])=\zeta_UD\simeq D/\mathfrak m$ as right
$D$-modules. This means that, for any $\alpha:U\to U[\infty]$, there
is a homomorphism $\beta\in D$ such that $\alpha=\zeta_U\beta$.
Moreover, if the above homomorphism $\alpha$ is non-zero, then
$\beta$ must be an isomorphism.

Keeping these details in mind, we can form the following commutative
diagram:
$$
 \xymatrix{(b)\quad &
U\ar[r]^-{\psi}\ar@{=}[d]&(R_\mathcal{U}/R)\otimes_RP_1\ar[d]_-{\simeq}\\
& U\ar[r]^-{(\zeta_U, g)}& U[\infty]\oplus E, }$$ where $E$ is an
$R$-module and $g:U\to E$ is an $R$-homomorphism which factorizes
through $\zeta_U$. Then, by applying Lemma \ref{lem3.6} and
combining $(a)$ with $ (b)$, we can construct the following exact
and commutative diagram:
$$\xymatrix{
0\ar[r]&U\ar[r]^-{\psi}\ar@{=}[d]&
(R_\mathcal{U}/R)\otimes_RP_1\ar[rr]^-{(R_\mathcal{U}/R)\otimes{_R}f_U}
\ar[d]^-{\simeq}&&(R_\mathcal{U}/R)\otimes_RP_0\ar[r]\ar[d]^-{\simeq}&0\\
0\ar[r]& U\ar[r]^-{(\zeta_U,\,0)}&U[\infty]\oplus
E\ar[rr]^-{\left({\pi_U\;\,0}\atop{\,0\;\;\;\,
1}\right)}&&(\tau^-U)[\infty]\oplus E\ar[r]&0 .}$$ Suppose
${\mathcal C} =\{U_1,U_2,\cdots, U_{m-1}, U_m\}$ with $m\geq 1$ such
that $\tau^-U_i=U_{i+1}$ for any $1\leq i\leq m$, where the
subscript of $U_i$ is always modulo $m$. Suppose $U=U_j$ for some
$1\leq j\leq m$. This means that $\pi_U$ coincides with
$\pi_{j,j+1}:U_j[\infty]\to U_{j+1}[\infty]$ defined in Lemma
\ref{lem3.3}(4), where $\pi_{m,m+1}=\pi_{m,1}$ by our convention.

Set $$M:=\bigoplus^m_{i=1}U_i[\infty],\quad \Lambda:=\End_R(M) \quad
\mbox{and} \quad \Pi:=\{\Hom_R(M,\pi_{s,s+1})\mid 1\leq s\leq m\}.$$

\textbf{Step $(2)$}. We prove $\Lambda_\Pi\simeq M_m\big(Q({\mathcal
C})\big)$, the $m\times m$ matrix ring over the division ring
$Q({\mathcal C})$.

For convenience, if $1\le u,v\le m$, we denote by $E_{u,v}$ the $
m\times m$ matrix unit which has $1$ in the $(u,v)$ position, and
$0$ elsewhere.

By Lemma \ref{lem3.4}, there is a ring isomorphism $\rho:
\;\Lambda\to\Gamma ({\mathcal C}),$ which sends
$E_{m,1}(\pi_{m,1})$ to $E_{m, 1}(\pi_{m,m})$ and
$E_{s,s+1}(\pi_{s,s+1})$ to $E_{s,s+1}$ for $1\leq s\leq m-1$ (see
Lemma \ref{lem3.4} for notations). Let $\varphi_m:
\Gamma(\mathcal{C})E_{m,m}\ra\Gamma(\mathcal{C})E_{1,1}$ and
$\varphi_s:
\Gamma(\mathcal{C})E_{s,s}\ra\Gamma(\mathcal{C})E_{s+1,s+1}$
be the canonical homomorphisms induced by multiplying on the right
by $E_{m, 1}(\pi_{m,m})$ and $E_{s,s+1}$ for $1\leq
s\leq m-1$, respectively, and define $\Theta:=\{\varphi_m\} \cup
\{\varphi_s\mid 1\leq s\leq m-1\}$. As a result, we get
$\Lambda_\Pi\simeq\Gamma(\mathcal{C})_\Theta$. It remains to
prove $\Gamma(\mathcal{C})_\Theta\simeq M_m\big(Q({\mathcal C})\big)$.
In fact, by Lemma \ref{lem2.3}, one can check that the canonical
inclusion from $\Gamma(\mathcal{C})$ to $M_m\big(D({\mathcal C})\big)$ is the
universal localization of $\Gamma(\mathcal{C})$ at
$\{\varphi_s\mid 1\leq s\leq m-1\}$. Observe that the universal
localization $D({\mathcal C})_{\pi_{m,m}}$ of $D({\mathcal C})$ at $\pi_{m,m}$ is equal to $Q(\mathcal C)$
by Lemma \ref{dv}. Now, combining Lemma \ref{lem2.5} with Corollary \cite[Corollary
3.4]{CX}, we have $$\Gamma(\mathcal{C})_\Theta\simeq
M_m\big(D({\mathcal C})\big)_{\varphi_m'}\simeq
M_m\big(\,D({\mathcal C})_{\pi_{m,m}}\,\big)\simeq M_m\big(Q(\mathcal C)\big),$$ where
$\varphi_m':M_m\big(D({\mathcal C})\big)E_{m,m}\to M_m\big(D({\mathcal C})\big)E_{1,1}$ is the
canonical homomorphism induced by $E_{m,1}(\pi_{m,m})$.  Thus
$\Lambda_\Pi\simeq M_m\big(Q(\mathcal C)\big)$.

\textbf{Step $(3)$}. We show that  $S_\Sigma$ is Morita equivalent
to the ad\`ele ring $\mathbb{A}_\mathcal{U}$ defined in Lemma
\ref{lem3.5}.

Indeed, by Lemma \ref{moddecom}(3), we have $R_\mathcal{U}/R\simeq
\bigoplus_{i\in I}\;\bigoplus_{V\in {\mathcal
C}_i}V[\infty]^{(\delta_V)},$ where
$\delta_{V}:=\dim_{\End_R(V)}\Ext^1_R(V,R)=
\dim_{\End_R(V)^{op}}(\tau V)\neq 0.$ We claim that there exists
$d\in\mathbb{N}$ such that $\delta_V\leq d $ for all $V\in
\mathcal{U}$.

In fact, let $\{S_j\mid 1\leq j\leq r\}$ be a complete set of
isomorphism classes of simple $R$-modules for some $r\in
\mathbb{N}$. For each $X\in R\modcat$, denote by
$\underline{\mbox{dim}}\,X\in\mathbb{N}^r$ the dimension vector of
$X$. Now, let $<-,->: \mathbb{N}^r\times\mathbb{N}^r\to \mathbb{Z}$
be the Euler form  of the tame hereditary $k$-algebra $R$, that is,
$<\underline{\mbox{dim}}\,Y,\,\underline{\mbox{dim}}\,Z>:=\dim_k\Hom_R(Y,Z)
-\dim_k\Ext^1_R(Y,Z)$ with $Y,Z\in R\modcat$, and further, let $q:
\mathbb{N}^r\to \mathbb{Z}$ be the quadratic form of $R$, that is,
$q(\underline{\mbox{dim}}\,Y):=<\underline{\mbox{dim}}\,Y,\,\underline{\mbox{dim}}\,Y>$,
and let $h=(h_i)_{1\leq i\leq r}$ be the minimal positive radical
vector of $q$. It is known that $h$ is equal to the sum of the
dimension vectors of all simple regular $R$-modules in ${\bf t'}$
for an arbitrary tube ${\bf t'}$ of $R$. Therefore, we have
$\delta_U\leq \dim_k (\tau U)\le
(\sum_{i}h_i)(\sum_{j}\dim_kS_j)<\infty$ for $U\in \mathscr{S}$. In
particular, if we take $d=(\sum_{i}h_i)(\sum_{j}\dim_kS_j)$, then
$\delta_V\le d$ for all $V\in \mathcal{U}$, as claimed.

Set
$$N:=\bigoplus_{i\in I}\;\bigoplus_{V\in
{\mathcal C}_i}V[\infty]\quad \mbox{and}\quad \Gamma:=\End_R(N).$$
By the above claim, one can check that $\Hom_R(R_{\mathcal U}/R,N)$
is a finitely generated, projective generator for $S\Modcat$, and
therefore $S$ is Morita equivalent to $\Gamma$. Note that Morita
equivalences preserve universal localizations by \cite[Corollary
3.4]{CX}.  Thus, we conclude from Step $(1)$ and the definition of
$\Sigma$ that $S_\Sigma$ is Morita equivalent to $\Gamma_\Phi$ with
$$\Phi:=\{\,\Hom_R(N,\pi_V)\mid V\in \mathcal{U}\}.$$

Now, let $\mathcal{U}=\mathcal{L}\dot{\cup} \mathcal{W}$ be an
arbitrary decomposition such that $\mathcal L$ is a union of cliques
${\mathcal C}_i$ with $i$ in an index set $I_0$ and that $\mathcal
W$ is a union of cliques ${\mathcal C}_j$ with $j$ in an index set
$I_1$. Note that $I=I_0\dot{\cup} I_1$. Moreover, if $i,j\in I$ with
$i\neq j$, then $\Hom_R(U[\infty],V[\infty])=0$ for all $U\in
{\mathcal C}_i$ and $V\in {\mathcal C}_j$. Thus, by Lemma
\ref{lem3.4}, we get the following isomorphisms:
$$(*)\quad \Gamma\simeq\prod_{i\in I}\End_R\big(\bigoplus_{V\in
{\mathcal C}_i}V[\infty]\big)\simeq\prod_{i \in I}\Gamma({\mathcal
C}_i)\simeq \prod_{i\in I_0}\Gamma({\mathcal C}_i) \times
\prod_{i\in I_1}\Gamma({\mathcal C}_i). $$ We write
$\Gamma_0:=\prod_{i\in I_0}\Gamma({\mathcal C}_i)$ and $
\Gamma_1:=\prod_{i\in I_1}\Gamma({\mathcal C}_i)$ and decompose
$\Phi=\Phi_0\cup\Phi_1$ where
$$\Phi_0:=\{\,\Hom_R(N,\pi_V)\mid V\in \mathcal{L}\,\}\quad\mbox{and}\quad
\Phi_1:=\{\,\Hom_R(N,\pi_W)\mid W\in\mathcal{W}\,\}.$$ Under these
isomorphisms $(*)$, we can regard $\Phi_0$ (respectively, $\Phi_1$)
as the set of homomorphisms between finitely generated projective
$\Gamma_0$-modules (respectively, $\Gamma_1$-modules). With these
identifications, one can prove
$\Gamma_\Phi\simeq(\Gamma_0)_{\Phi_0}\times (\Gamma_1)_{\Phi_1}.$

Next, we assume that  each clique in $\mathcal{W}$ is of rank one,
and each clique $L\in\mathcal{L}$ is of rank greater than one.
Clearly, $\mathcal L$ is a finite set.

On the one hand, by the foregoing discussion and Step $(2)$, we
obtain
$$(\Gamma_0)_{\Phi_0}\simeq\prod_{i\in I_0}M_{c(\mathcal{C}_i)}\big(\,Q({\mathcal{C}_i})\,\big).$$

On the other hand, we have $\Gamma_1=\prod_{i\in
I_1}\,D({\mathcal{C}_i})\,$. Now, we claim
$(\Gamma_1)_{\Phi_1}\simeq\mathbb{A}_{\,\mathcal{W}}$, where
$$\mathbb{A}_{\,\mathcal{W}}:=\bigg\{(f_{i})_{i\in I_1}\in\prod_{i\in I_1}Q({\mathcal{C}_i})
\,\mid f_{i}\in D({\mathcal{C}_i})\; \mbox{for almost all }
i\in{I_1}\bigg\}.$$ This ring is  similar to the so called ad\`ele
ring appearing in the algebraic number theory (see \cite[Chapter 5,
Section 1]{ne}).

Actually, for each $i\in I_1$,  the clique $\mathcal{C}_i$ consists
of only one simple regular module. Hence we write
$D({\mathcal{C}_i})=\End_R(\mathcal{C}_i)$, which is a discrete
valuation ring with a unique maximal ideal generated by $\pi_{i}$.

We define $e_{i}:=\big(\beta_{j}\big)_{j\in I_1}\in\Gamma_1$ by
$\beta_{i}=1$ and $\beta_{j}=0$ if $j\neq i$, and define
$\varepsilon_{i}:=\big(\theta_{j}\big)_{j\in I_1}\in\Gamma_1$ by
$\theta_{i}=\pi_i$ and $\theta_{j}=1$ if $j\neq i$. Let
$\varphi_{i}:\Gamma_1 e_{i}\to\Gamma_1e_{i}$ be the right
multiplication map defined by $g\mapsto g\pi_i$ for any $g\in
D({\mathcal{C}_i})$.

Under the isomorphisms $(*)$, we can identify $\Phi_1$ with
$\{\varphi_{j}\mid j\in I_1\}$. Note that the right multiplication
map $\ol{\varepsilon}_i$ defined by $\varepsilon_i$ has the
following form:
$$\ol{\varepsilon}_i=\left({\varphi_i\;\,0}\atop{\,0\;\;\;\, 1}\right):\;
\Gamma_1e_i\oplus \Gamma_1(1-e_i)\lra \Gamma_1e_i\oplus
\Gamma_1(1-e_i).$$ Set  $\Psi:=\{\ol{\varepsilon}_j\mid j\in I_1\}$.
It is easy to see that $(\Gamma_1)_{\Phi_1}$ is isomorphic to the
universal localization $(\Gamma_1)_{\,\Psi}$ of $\Gamma_1$ at
$\Psi$. We consider the minimal multiplicative subset $\Upsilon$ of
$\Gamma_1$ containing all $\varepsilon_j$ for $j\in I_1$. Clearly,
$(\Gamma_1)_{\Psi}$ is isomorphic to the universal localization of
$\Gamma_1$ at $\Upsilon$, that is, the universal localization of
$\Gamma_1$ at all  right multiplication maps induced by the elements
of $\Upsilon$.  One can check
$$\Upsilon=\bigg\{(f_{i})_{i\in
I_1}\in\prod_{i\in I_1}\big\{(\pi_i)^n\mid n\in{\mathbb N}\big\}
\;\bigg|\; f_{i}=1\; \mbox{for almost all }
i\in{I_1}\bigg\}\subseteq \Gamma_1.$$ We claim that  $\Upsilon$ is a
left and right denominator subset of $\Gamma_1$ (see Definition
\ref{denom}).

Indeed, let $a=(a_i)_{i\in I_1}\in\Gamma_1$ and
$s=(\pi_i^{n_i})_{i\in I_1}\in\Upsilon$ with $n_i\in\mathbb{N}$.
Since $D({\mathcal{C}_i})$ is a discrete valuation ring for each
$i\in I_1$,  we have
$D({\mathcal{C}_i})\pi_i^{n_i}=\pi_i^{n_i}D({\mathcal{C}_i})$, and
therefore $\Gamma_1s=\prod_{i\in
I_1}D({\mathcal{C}_i})\pi_i^{n_i}=\prod_{i\in
I_1}\pi_i^{n_i}D({\mathcal{C}_i})$. This means  $sa\in\Upsilon a\cap
\Gamma_1s\neq\emptyset$, which verifies the condition $(i)$ in
Definition \ref{denom}. On the other hand, if $as=0$, then
$a_i\pi_i^{n_i}=0$. Since $\pi_i^{n_i}\neq 0$ and
$D({\mathcal{C}_i})$ is a domain for $i\in I_1$,  we have $a_i=0$,
and so $a=0$, which verifies the condition $(ii)$ in Definition
\ref{denom}. Thus, $\Upsilon$ is a left denominator subset of
$\Gamma_1$. Similarly, we can prove that $\Upsilon$ is also a right
denominator subset of $\Gamma_1$.

It remains to prove
$\Upsilon^{-1}\Gamma_1\simeq\mathbb{A}_{\mathcal{W}}$. In fact, it
follows from Lemma \ref{Ore} that the universal localization of
$\Gamma_1$ at $\Upsilon$ is the same as the Ore localization
$\Upsilon^{-1}\Gamma_1$ of $\Gamma_1$ at $\Upsilon$. Moreover, by
Lemma \ref{dv}, we see that, for each $j\in I_1$, the Ore
localization of $D({\mathcal{C}_j})$ at $\{(\pi_j)^n\mid
n\in{\mathbb N}\}$ is the division ring $Q({\mathcal{C}_j})$ of
fractions of $D({\mathcal{C}_j})$. Thus, by the definition of Ore
localizations (see Lemma \ref{Ore}), one can easily prove
$\Upsilon^{-1}\Gamma_1\simeq\mathbb{A}_{\mathcal{W}}$.

Summing up what we have proved, we get
$$\Gamma_\Phi\simeq(\Gamma_0)_{\Phi_0}\times
(\Gamma_1)_{\Phi_1}\simeq\prod_{i\in I_0}\,M_{c({\mathcal
C}_i)}\big( Q({\mathcal{C}_i})\big) \; \times \;
\mathbb{A}_{\mathcal W}, $$ the latter is Morita equivalent to
$\mathbb{A}_{\mathcal U}$.  As $S_\Sigma$ is Morita equivalent to
$\Gamma_\Phi$, we see that $S_\Sigma$ is Morita equivalent to
$\mathbb{A}_\mathcal{U}$. This completes the whole proof. $\square$

\bigskip
{\bf Proof of Theorem \ref{th1.1}.} Recall that
$B=\End_R(R_{\mathcal U}\oplus R_{\mathcal U}/R)$ and
$S:=\End_R(R_{\mathcal U}/R)$. By Corollary \ref{lem2.4}, there is a
recollement of derived module categories:

$$(*)\quad \xymatrix@C=1.2cm{\D{S_{\Sigma}}\ar[r]&\D{B}\ar[r]
\ar@/^1.2pc/[l]\ar@/_1.2pc/[l]
&\D{R}\ar@/^1.2pc/[l]\ar@/_1.2pc/[l]},\vspace{0.3cm}$$ where
$S_\Sigma$ is the universal localization  of $S$ at
$\Sigma:=\{S\otimes_Rf_U\mid U\in {\mathcal U}\}$.

Now we write $\mathcal{U}=\mathcal{U}_0\cup
\mathcal{U}_1\subseteq{\mathscr S}$ such that $\mathcal{U}_0$
contains no cliques and $\mathcal{U}_1$ is a union of cliques
${\mathcal C}_i$ with $i\in I$, an index set. We conclude from Lemma
\ref{lem3.2} that $S_\Sigma$ is isomorphic to the universal
localization $\Lambda_\Theta$ of $\Lambda$ at $\Theta$ with
$\Lambda:=\End_{R_{\mathcal{U}_0}}(R_\mathcal{U}/R_{\mathcal{U}_0})$
and
$\Theta:=\{\Lambda\otimes_{R_{\mathcal{U}_0}}({R_{\mathcal{U}_0}}\otimes_Rf_V)\mid
V\in \mathcal{U}_1\}$. Note that $R_{\mathcal{U}_0}$ is a
finite-dimensional tame hereditary $k$-algebra, and that
$\mathcal{U}_1$ is a union of cliques when regarded as a set of
simple regular $R_{\mathcal{U}_0}$-modules. Now, by applying Lemma
\ref{lem3.5} to $R_{\mathcal{U}_0}$ and $\mathcal{U}_1$, we can
deduce that $\Lambda_\Theta$ is Morita equivalent to the ad\`ele
ring $\mathbb{A}_{\mathcal U}$ in Theorem \ref{th1.1}.

Thus, we have proved that $S_\Sigma$ is Morita equivalent to
$\mathbb{A}_{\mathcal U}.$  If we substitute $\D{S_\Sigma}$ by
$\D{\mathbb{A}_{\mathcal U}}$ in $(*)$, then we  obtain the desired
recollement of derived module categories in Theorem \ref{th1.1}:

$$\xymatrix@C=1.2cm{\D{\mathbb{A}_{\mathcal U}}\ar[r]&\D{B}\ar[r]
\ar@/^1.2pc/[l]\ar@/_1.2pc/[l]
&\D{R}\ar@/^1.2pc/[l]\ar@/_1.2pc/[l]}.\vspace{0.3cm}$$ This
completes the proof of the first part of Theorem \ref{th1.1}.

As for the second part, we note that, if $k$ is algebraically
closed, then, for each clique $\mathcal{C}$ of $R$, the rings
$D({\mathcal{C}})$ and $Q({\mathcal{C}})$ are isomorphic to $k[[x]]$
and $k((x))$ by Lemma \ref{lem3.3}(5), respectively. Now, combining
this with the first part of Theorem \ref{th1.1}, we know that
$\mathbb{A}_{\mathcal U}$ is isomorphic to $\mathbb{A}_I$. This
finishes the  proof. $\square$

\bigskip
If we take $\mathcal U$ = $\mathscr S$, then the tilting $R$-module
$R_{\mathscr S}\oplus R_{\mathscr S}/R$ is a Reiten-Ringel tilting
module (see \cite{R}).  This tilting module is actually of the form
$G^{(n)}\oplus \bigoplus_{U\in {\mathscr S}}U[\infty]^{(\delta_U)}$,
where $G$ is the unique generic $R$-module with
$n=\dim_{\End_R(G)}G$, and
$\delta_{U}=\dim_{\End_R(U)}\Ext^1_R(U,R)$ for $U\in \mathscr{S}$
(see \cite[Proposition 1.8]{HJ2}). Recall that $\mathscr S$ is
parameterized by the projective line ${\mathbb P}^1(k)$ if $k$ is
algebraically closed.  As a consequence of Theorem \ref{th1.1}, we
have the following corollary.

\begin{Koro}\label{Adele}
If $k$ is an algebraically closed field and $T$ is the Reiten-Ringel
tilting $R$-module $T_{\mathscr S}$, then there is a recollement

$$\xymatrix@C=1.2cm{\D{\mathbb{A}_{\mathbb{P}^1(k)}}\ar[r]&\D{\End_R(T)}\ar[r]
\ar@/^1.2pc/[l]\ar@/_1.2pc/[l]
&\D{R}\ar@/^1.2pc/[l]\ar@/_1.2pc/[l]}.$$
\end{Koro}

\section{Stratifications of derived module categories}

In this section, we shall use Theorem \ref{th1.1} to get
stratifications of the derived categories of the endomorphism rings
of tilting modules of the form $R_{\mathcal U}\oplus R_{\mathcal
U}/R$. It turns out that our consideration for general tame
hereditary algebras  is converted into understanding the case of
special tame hereditary algebras consisting of two isomorphism
classes of simple modules. In particular, if $k$ is an algebraically
closed field,  we are led to the Kronecker algebra. In this way, we
shall prove Corollary \ref{cor1.2} in this section.

\subsection{Universal localizations of general tame hereditary algebras}

In this subsection, we shall discuss the endomorphism algebras of
tilting modules  associated with universal localizations of tame
hereditary algebras at simple regular modules. The consideration
here will be served as a part of preparations for stratifications of
derived categories in Subsection 4.3.

Throughout this subsection,  $R$ is an  indecomposable
finite-dimensional tame hereditary algebra over an arbitrary field
$k$, and ${\mathscr S}:={\mathscr S}(R)$ is the complete set of
isomorphism classes of all simple regular $R$-modules.

Let $\mathcal{U}$ be an arbitrary subset of  $\mathscr{S}$. The
following result gives a characterization of the universal
localization $R_{\mathcal U}$ of $R$ at $\mathcal U$ from  the view
of derived equivalences.

\begin{Lem}
Let $\mathcal{U}\subseteq{\mathscr S}$. Then there exists
$\mathcal{V}\subseteq{\mathscr S}$ with
$\mathcal{U}\cap\mathcal{V}=\emptyset$ such that, for ${\mathcal
W}:=\mathcal{U}\cup\mathcal{V}$, the following statements are true.

$(1)$ There is a finite-dimensional tame hereditary $k$-algebra
$\Lambda$ with only two non-isomorphic simple modules, and a set
$\mathcal{S}$ of simple regular $\Lambda$-modules such that
$R_{{\mathcal W}}$ coincides with the universal localization
$\Lambda_{\mathcal{S}}$ of $\Lambda$ at $\mathcal{S}$.

$(2)$ The $R_\mathcal{U}$-module $T:=R_{{\mathcal W}}\oplus
R_{{\mathcal W}}/R_\mathcal{U}$ is a classical tilting module. In
particular, $R_{\mathcal{U}}$ and $\End_{R_{\mathcal U}}(T)$ are
derived-equivalent.\label{lem4.2}
\end{Lem}

{\it Proof.} Suppose $\mathcal{U}=\mathcal{U}_0\dot{\cup}
\mathcal{U}_1\subseteq{\mathscr S}$ such that $\mathcal{U}_0$
contains no cliques and $\mathcal{U}_1$ is a union of cliques.
Observe that we may assume $\mathcal{U}_0=\emptyset.$ In fact, if
$\mathcal{U}_0$ is not empty, we can replace $R$ by
$R_{\mathcal{U}_0}$ and $\mathcal{U}$ by $\mathcal{U}_1$ since
$R_{\mathcal{U}_0}$ is a tame hereditary algebra and $\mathcal{U}_1$
can be seen as a set of simple regular $R_{\mathcal{U}_0}$-modules.

From now on, we suppose $\mathcal{U}_0=\emptyset,$  that is,
$\mathcal{U}$ is a union of cliques.  Let $\mathcal{V}$ be a maximal
subset of ${\mathscr S}$  with respect to the following property:
$\mathcal{V}\cap\mathcal{U}=\emptyset$ and $\mathcal{V}$ contains no
cliques. In other words, from each clique $\mathcal C$ not contained
in $\mathcal{U}$, we choose $c({\mathcal C})-1$ elements, and let
$\mathcal{V}$ be the union of all these elements. Clearly, the
choice of $\mathcal{V}$ is not unique in general.

Let ${\mathcal W}:=\mathcal{U}\,\dot{\cup}\mathcal{V}$, and let
${\mathcal U}_{>1}$ be the union of all cliques ${\mathcal C}_{i\in
I}$ in $\mathcal U$ of rank greater than one, where $I$ is a finite
set. We choose $c(\mathcal{C}_i)-1$ elements from each
$\mathcal{C}_i$ for $i\in I$, and let ${\mathcal V}'$ be the set
consisting of all of these elements. Now, we define
$\mathcal{L}={\mathcal V}\cup {\mathcal V}'$ and write ${\mathcal
W}=\mathcal{L}\dot{\cup}\mathcal{M}$.

We claim that the statement $(1)$ holds true. Indeed, it follows
from Lemma \ref{lem3.1}(1) that $R_\mathcal{L}$ is a tame hereditary
algebra such that all cliques of $R_\mathcal{L}$ consist of only one
simple regular module. This means that $R_{\mathcal{L}}$ has exactly
two isomorphism classes of simple modules.  By Lemma
\ref{lem3.1}(3), we have $R_{{\mathcal
W}}=(R_\mathcal{L})_{\ol{\mathcal{M}}}$ with
$\ol{\mathcal{M}}:=\{R_\mathcal{L}\otimes_R L\mid
L\in{\mathcal{M}}\}$. Thus, setting $\Lambda:=R_\mathcal{L}$ and
$\mathcal{S}:={\ol{\mathcal{M}}}$, we get the statement $(1)$.

In the following, we shall show the statement $(2)$. Note that
$\mathcal{V}$ contains no cliques. Thus, it follows from Lemma
\ref{lem3.1}(1) that $R_\mathcal{V}$ is a tame hereditary algebra
and $R_\mathcal{V}/R$ is a finitely presented $R$-module.  By Lemma
\ref{3.4}(1), $R_{{\mathcal W}}/R_\mathcal{U}\simeq
R_{\mathcal{U}}\otimes_R (R_{\mathcal{V}}/R)$ as
$R_\mathcal{U}$-$R$-bimodules. This implies that $R_{{\mathcal
W}}/R_\mathcal{U}$ is a finitely presented $R_\mathcal{U}$-module,
and so is the $R_\mathcal{U}$-module $T$. Hence, $T$ is a classical
$R_\mathcal{U}$-module.  $\square$

\medskip
As a consequence of Lemma \ref{lem4.2}, we obtain the following
result, which describes  $R_{\mathcal{U}}$ up to derived equivalence
by a triangular matrix ring such that the rings in the diagonal are
relatively simple.

\begin{Koro}
Suppose that  $\mathcal{U}\subseteq{\mathscr S}$ is a union of
cliques ${\mathcal C}_{i\in I}$ with $I$ an index set. Let
$\mathcal{V}$ be a maximal subset of ${\mathscr S}$ such that
$\mathcal{V}\cap\mathcal{U}=\emptyset$ and $\mathcal{V}$ contains no
cliques, and let ${\mathscr C}({\mathcal V})=\dot{\cup}_{j\in
J}{\mathcal C}_j$ with $J$ an index set.  Define ${\mathcal
W}:=\mathcal{U}\cup\mathcal{V}$ and
$T_\mathcal{U}:=R_{\mathcal{U}}\oplus R_{\mathcal{U}}/R$.  Then the
following statements hold true:

$(1)$ There is a canonical ring isomorphism:
$$\End_R(T_\mathcal{U})\simeq\left(\begin{array}{lc} R_\mathcal{U} &
\Hom_R(R_\mathcal{U},\,
R_\mathcal{U}/R)\\
\,0&\End_R(R_\mathcal{U}/R)\end{array}\right).$$

$(2)$  $\End_R(R_\mathcal{U}/R)$ is Morita equivalent to
$\prod_{i\in I}\Gamma\big(\mathcal{C}_i\big)$, where
$\Gamma(\mathcal{C})$ is defined in Lemma \ref{lem3.4} for each
clique $\mathcal{C}$ of $R$.

$(3)$ $R_\mathcal{U}$ is derived-equivalent to the following
triangular matrix ring
$$
\End_{R_{\mathcal U}}(R_{{\mathcal W}}\oplus R_{{\mathcal
W}}/R_\mathcal{U})=\left(\begin{array}{lc} R_{{\mathcal W}} &
\Hom_{R_\mathcal{U}}\big(R_{{\mathcal W}},\,
R_{{\mathcal W}}/R_\mathcal{U}\big)\\
\,0&\End_{R_\mathcal{U}}\big(R_{{\mathcal
W}}/R_\mathcal{U}\big)\end{array}\right)
$$
such that

$(a)$ $R_{{\mathcal W}}$ is the universal localization
$\Lambda_{\mathcal{S}}$ of a finite-dimensional tame hereditary
$k$-algebra $\Lambda$, which has two isomorphism classes of simple
modules, at a set $\mathcal{S}$ of simple regular $\Lambda$-modules,
and

$(b)$ $\End_{R_\mathcal{U}}(R_{{\mathcal W}}/R_\mathcal{U})$ is
Morita equivalent to $\prod_{j\in J}T_{c({\mathcal
C}_j)-1}\big(\,\End_R(V_j)\,\big)$, where $V_j\in\mathcal{C}_j$ is a
fixed element for each $j\in J$, and $T_n(A)$ stands for the
$n\times n$ upper triangular matrix ring over a ring
$A$.\label{cor4.3}
\end{Koro}

{\it Proof.} Clearly, $(1)$ follows from $\lambda:R\to
R_\mathcal{U}$ being a ring epimorphism and
$\Hom_R(R_\mathcal{U}/R,R_\mathcal{U})=0$. $(2)$ follows from ($*$)
in Step $(3)$ of the proof of Lemma \ref{lem3.5}. As to $(3)$, we
first show the statement $(b)$. In fact, by the proof of Lemma
\ref{lem4.2}, we know $R_{{\mathcal W}}/R_\mathcal{U}\simeq
R_{\mathcal{U}}\otimes_R (R_{\mathcal{V}}/R)$ as
$R_\mathcal{U}$-$R$-bimodules. Since $\mathcal{V}\subseteq
\mathcal{U}^{+}$, we have $R_{\mathcal{U}}\otimes_R
(R_{\mathcal{V}}/R)\simeq R_{\mathcal{V}}/R$ as $R$-modules by Lemma
\ref{orth}, and therefore $R_{{\mathcal W}}/R_\mathcal{U}\simeq
R_{\mathcal{V}}/R$ as $R$-modules. This implies that
$\End_{R_\mathcal{U}}(R_{{\mathcal W}}/R_\mathcal{U})\simeq
\End_R(R_{{\mathcal W}}/R_\mathcal{U})\simeq
\End_R(R_\mathcal{V}/R).$ Now, by Lemma \ref{moddecom}(2), one can
prove
$$R_{\mathcal{V}}/R\simeq
\bigoplus_{j\in J}\bigoplus_{i=1}^{c({\mathcal C}_j)-1}U_{i,j}
\,[c({\mathcal C}_j)-i]^{(\delta_{{i,j}})},$$ where
$\delta_{{i,j}}>0$ and $\mathcal{V}\cap {\mathcal C}_j=\{U_{i,j}\mid
1\leq i< c({\mathcal C}_j)\}$ such that $U_{i+1,\, j}=\tau^-U_{i,j}$
for all $1\leq i<c({\mathcal C}_j)-1.$ Further, for any $j\in J$,
one can check
$$\End_R\big(\bigoplus_{i=1}^{c({\mathcal C}_j)-1}U_{i,j}
\,[c({\mathcal C}_j)-i]\;\big)\simeq T_{c({\mathcal
C}_j)-1}\big(\,\End_R(V_j)\,\big),$$ where $V_j$ is a fixed element
of $\mathcal{C}_j$ with $j\in J$. Note that $\End_R(V_j)$ is
independent of the choice of elements of $\mathcal{C}_j$ up to
isomorphism. Thus $\End_{R_\mathcal{U}}(R_{{\mathcal
W}}/R_\mathcal{U})$ is Morita equivalent to $\prod_{j\in
J}T_{c({\mathcal C}_j)-1}\big(\End_R(V_j)\big)$, since there is no
non-trivial homomorphism between two different tubes.

Note that the other conclusions in $(3)$ are  consequences of Lemma
\ref{lem4.2} and of properties of injective ring epimorphisms (see
also \cite[Lemma 6.4(2)]{CX}). This completes the proof. $\square$

\bigskip
Thus, by Corollary \ref{cor4.3}(3), the consideration of the derived
category $\D{R_{\mathcal U}}$ needs first to understand universal
localizations of  tame hereditary algebras with two isomorphism
classes of simple modules, at simple regular modules. If $k$ is an
algebraically closed field, then each tame hereditary algebra with
two isomorphism classes of simple modules is Morita equivalent to
the Kronecker algebra. So, in the next subsection, we shall focus
our attention on the universal localizations of the Kronecker
algebra.

\subsection{Universal localizations of the Kronecker algebra at simple regular modules}

In this subsection,  we shall consider the particular tame
hereditary algebra, the Kronecker algebra.  The  results obtained
here will be served again as a preparation for the discussion of
stratifications of derived module categories in the next subsection.

Throughout this subsection, $k$ is a field, and $R$ is the Kronecker
algebra
$\left(\begin{array}{cc} k & k^2\\
0 & k\end{array}\right)$, where the $k$-$k$-bimodule structure of
$k^2$ is given by $a(b, c)d=(abd, acd)$ with $a, b, c, d\in k$.
 It is known that $R$
can be interpreted as the path algebra of the quiver
$$Q:\;
\xymatrix{2\ar@<0.4ex>[r]^-{\alpha}\ar@<-0.4ex>[r]_-{\beta}& 1},$$
and that $R\Modcat$ (respectively, $R\modcat$) is equivalent to the
category of (respectively, finite-dimensional) representations of
$Q$ over $k$.

In this subsection, we denote by $V$ the representation
$\xymatrix{k\ar@<0.4ex>[r]^-{0}\ar@<-0.4ex>[r]_-{1}& k}$.  By Lemma
\ref{lem2.3}, one can check that $R_{V}=M_2(k[x])$, and the
universal localization $\lambda: R\to R_{V}$ is given by
$\left(\begin{array}{cc} a& (c,d)\\
0 & b\end{array}\right)\mapsto\left(\begin{array}{cc} a
& c+dx\\
0 & b \end{array}\right)$ for $a,b,c, d\in k$. In particular, the
restriction functor $\lambda_*:R_V\Modcat\to R\Modcat$ induced by
$\lambda$ is fully faithful. Let $e=\left(\begin{array}{ll} 1 &0\\
0 & 0\end{array}\right)\in R_V$. Clearly, the tensor  functor
${R_V}\,e\otimes_{k[x]}-: \,k[x]\Modcat\to R_V\Modcat$ is an
equivalence. Now, we define $F: k[x]\Modcat\to R\Modcat$ to be the
composition of the functors ${R_V}\,e\otimes_{k[x]}-$ and
$\lambda_*$. Then $F$ is a fully faithful exact functor, and sends
each $k[x]$-module $M$ to the representation
$\xymatrix{M\ar@<0.4ex>[r]^-{1}\ar@<-0.4ex>[r]_-{x}& M}$. Moreover,
we have the following result.

\begin{Lem}{\rm \cite[Theorem 4]{R0}} \label{R0}
The functor $F$ induces an equivalence between the category of
finite-dimensional $k[x]$-modules and  the category of
finite-dimensional regular $R$-modules with regular composition
factors not isomorphic to $V$.
\end{Lem}

%By Lemma \ref{R0}, there is a bijection between the set of
%isomorphism classes of simple $k[x]$-modules and that of simple
%regular $R$-modules which are not isomorphic to $V$. Note that each
%simple $k[x]$-module corresponds to a monic irreducible polynomial
%in $k[x]$.

Let $\mathcal{P}$ be the set of all monic irreducible polynomials in
$k[x]$. For each $p(x)\in\mathcal{P}$, we denote by $k_{p(x)}$ the
extension field $k[x]/(p(x))$ of $k$, and by $V_{p(x)}$ the
representation
$\xymatrix{k_{p(x)}\ar@<0.4ex>[r]^-{1}\ar@<-0.4ex>[r]_-{x}&
k_{p(x)}}$, which is the image of $k_{p(x)}$ under $F$. Since simple
$k[x]$-modules are parameterized by monic irreducible polynomials,
it follows from Lemma \ref{R0} that ${\mathscr
S}:=\{V\}\cup\{V_{p(x)}\mid p(x)\in\mathcal{P}\}$ is a complete set
of isomorphism classes of simple regular $R$-modules. If $k$ is
algebraically closed, then $\mathcal{P}=\{x-a\mid a\in k\}$, and
therefore $\mathscr S$ can be identified with the projective line
${\mathbb P}^1(k)$.

\medskip
The following corollary gives a characterization of the
endomorphisms rings of Pr\"ufer modules.

\begin{Koro}
Let $t$ be a variable and $p(x)\in\mathcal{P}$. Then there are
isomorphisms of rings:
$$\End_R\big(V[\infty]\big)\simeq k[[t]] \quad \mbox{and}\quad
\End_R\big(V_{p(x)}[\infty]\big)\simeq k_{p(x)}\,[[t]].$$
\label{prufer}
\end{Koro}

{\it Proof.} Recall that, for any simple regular $R$-module $U$, we
have $\End_R\big(U[\infty]\big)\simeq
\displaystyle\varprojlim_n\,\End_R(U[n])$ as rings. If $U=V$, then
$\End_R(U[n])\simeq k[t]/(t^n)$ for any  $n>0$, and therefore
$\End_R\big(U[\infty]\big)\simeq\displaystyle\varprojlim_n\,k[t]/(t^n)\simeq
k[[t]]$. Suppose  $U=V_{p(x)}$. It follows from Lemma \ref{R0} that
$U[n]\simeq F\big(\,k[x]/(p(x)^n)\,\big)$ as $R$-modules, and that
$\End_R(U[n])\simeq\End_{k[x]}\big(\,k[x]/(p(x)^n)\,\big)\simeq
k[x]/(p(x)^n)$ for any $n>0$. Thus $\End_R\big(U[\infty]\big)\simeq
\displaystyle\varprojlim_n\,k[x]/(p(x)^n)$. This implies that
$\End_R\big(U[\infty]\big)$ is a complete commutative discrete
valuation ring (see Lemma \ref{lem3.3}(5)), and therefore it is a
regular ring of Krull dimension 1. Recall that a regular ring is by
definition a commutative noetherian ring of finite global dimension.
For regular rings, the global dimension agrees with the Krull
dimension.

It remains to prove $\displaystyle\varprojlim_n\,k[x]/(p(x)^n)\simeq
k_{p(x)}\,[[t]]$. Actually, this follows straightforward from the
following classical result (see \cite[Theorem 15]{Cohen} for
details):

Let $S$ be a complete regular local ring of Krull dimension $m$ with
the residue class field $K$. If $S$ contains a field, then $S$ is
isomorphic to  the formal power series ring $K[[t_1,\cdots,t_m]]$
over $K$ in variables $t_1,\cdots,t_m$.

\medskip
Hence $\End_R\big(U[\infty]\big)\simeq
\displaystyle\varprojlim_n\,k[x]/(p(x)^n)\simeq k_{p(x)}\,[[t]]$,
which finishes the proof.  $\square$

\smallskip
In the remainder of this subsection, let $\Delta$ be a subset of
$\mathcal{P}$, and let $\mathcal{U}:=\{V\}\cup\{V_{p(x)}\mid p(x)\in
\Delta\}$. We define the $\Delta$-ad\`ele ring of $k[x]$ as follows:
$$\mathbb{A}(\Delta):=k((t)) \times \bigg\{\big(\theta_{p(x)}\big)_{p(x)\in
\Delta}\, \in\prod_ {p(x)\in
\Delta}k_{p(x)}\,((t))\;\big|\;\theta_{p(x)}\in k_{p(x)}\,[[t]]\;
\mbox{ for almost all }\; p(x)\in \Delta\bigg\}.$$

Combining Theorem \ref{th1.1} with Corollary \ref{prufer}, we get
the following result.

\begin{Koro}\label{adele ring}
Let $B$ be the endomorphism ring of the tilting $R$-module
$R_{\mathcal{U}}\oplus R_{\mathcal{U}}/R$.  Then there is a
recollement of derived categories:

$$
\xymatrix@C=1.2cm{\D{\mathbb{A}({\Delta})}\ar[r]&\D{B}\ar[r]
\ar@/^1.2pc/[l]\ar@/_1.2pc/[l]
&\D{R}\ar@/^1.2pc/[l]\ar@/_1.2pc/[l]}.\vspace{0.3cm}
$$
\end{Koro}

\medskip
In Corollary \ref{adele ring}, if $\Delta=\mathcal{P}$, then the
$\mathcal{P}$-ad\`ele ring $\mathbb{A}({\mathcal{P}})$ of $k[x]$
coincides with the ad\`ele ring $\mathbb{A}_{k(x)}$ of the fraction
field $k(x)$, which appears in global class field theory (see
\cite[Chapter VI]{ne} and \cite[Theorem 2.1.4]{ep}).

\medskip
Finally, we prove the following lemma as the last preparation for
the proof of Corollary \ref{cor1.2}.

\begin{Lem} Let $D$
be the smallest subring of the fraction field $k(x)$ of $k[x]$
containing both $k[x]$ and $\frac{1}{p(x)}$ with  all
$p(x)\in\Delta$. Then $R_{\mathcal{U}}\simeq M_2(D)$, the $2\times
2$ matrix ring over $D$. In particular, $R_{\mathcal{U}}$ is Morita
equivalent to the Dedekind integral domain $D$. \label{lem4.1}
\end{Lem}

{\it Proof.} Define $\mathcal{W}:=\{R_V\otimes_R V_{p(x)} \mid
p(x)\in\Delta\}$. Then $R_{\mathcal{U}}=(R_V)_{\mathcal{W}}$ by
Lemma \ref{lem3.1}(3). Recall that $R_V=M_2(k[x])$ and $\lambda:
R\to R_{V}$ is the universal localization of $R$ at $V$. On the one
hand, for each $p(x)\in\Delta$, it follows from
$V_{p(x)}=F\big(k_{p(x)}\big)=
\lambda_*\big({R_V}\,e\otimes_{k[x]}k_{p(x)}\big)$ that
$$R_V\otimes_RV_{p(x)}\simeq V_{p(x)}={R_V}\,e\otimes_{k[x]}k_{p(x)}
=\left(\begin{array}{c} k_{p(x)}\\ k_{p(x)}\end{array}\right)$$ as
$R_V$-modules. On the other hand, by \cite[Corollary 3.4]{CX},
Morita equivalences preserve universal localizations. Consequently,
we have $R_\mathcal{U}=\big(M_2(k[x])\big)_{\mathcal W}\,\simeq\,
M_2\big(k[x]_{\Theta}\big)$ with $\Theta:=\{k_{p(x)}\mid
p(x)\in\Delta\}\subseteq k[x]\Modcat.$ Now, one may readily see that
$k[x]_{\Theta}$ coincides with the localization of $k[x]$ at the
smallest multiplicative subset of $k[x]$ containing $\{p(x)\mid
p(x)\in\Delta\}$, which is exactly the ring $D$ defined in Lemma
\ref{lem4.1}. Since $k[x]$ is a Dedekind integral  domain and since
localizations of Dedekind  integral domains are again Dedekind
integral domains, we see that $D$ is a Dedekind integral domain. As
a result, we have $R_\mathcal{U}\simeq M_2(D)$. This completes the
proof. $\square$

\medskip
{\it Remarks.}  (1) If $k$ is an algebraically closed field, then,
for any simple regular $R$-module $U$, we can choose an automorphism
$\sigma:R\to R$, such that the induced functor $\sigma_*:
R\Modcat\to R\Modcat$ by $\sigma$ is an equivalence with
$\sigma_*(U)\simeq V$. This implies that, up to isomorphism,  Lemma
\ref{lem4.1} provides a complete description of $R_{\mathcal{V}}$
for any subset $\mathcal{V}$ of $\mathscr{S}$. In particular,
$R_{\mathcal{V}}$  is Morita equivalent to a Dedekind integral
domain.

(2) If we localize $R$ at all non-isomorphic simple regular modules
${\mathscr S}$ which is indexed by all monic irreducible
polynomials, then, by Lemma \ref{lem4.1}, we have $R_{\mathscr
S}\simeq M_2(k(x))$ since the smallest subring containing the
inverses of all irreducible polynomials $p(x)$ is just $k(x)$.

\subsection{Stratifications of derived module categories}

The main purpose of this subsection is to prove Corollary
\ref{cor1.2}. We first recall the definition of stratifications of
derived categories of rings.

As in \cite{HKL}, the derived module category $\D{A}$ of a ring $A$
is called derived simple if it is not a non-trivial recollement of
any derived categories of rings. A stratification of $\D{A}$ of a
ring $A$ by derived categories of rings is defined to be a sequence
of iterated recollements of the following form: a recollement of
$A$, if it is not derived simple,
$$\xymatrix@C=1.2cm{\D{A_1}\ar[r]&\D{A}\ar[r]
\ar@/^1.2pc/[l]\ar@/_1.2pc/[l]
&\D{A_2}\ar@/^1.2pc/[l]\ar@/_1.2pc/[l]},\vspace{0.3cm}$$  a
recollement of the ring $A_1$, if it is not derived simple,
$$\xymatrix@C=1.2cm{\D{A_{11}}\ar[r]&\D{A_1}\ar[r]
\ar@/^1.2pc/[l]\ar@/_1.2pc/[l]
&\D{A_{12}}\ar@/^1.2pc/[l]\ar@/_1.2pc/[l]},\vspace{0.3cm}$$ and a
recollement of the ring $A_2$, if it is not derived simple,
$$\xymatrix@C=1.2cm{\D{A_{21}}\ar[r]&\D{A_2}\ar[r]
\ar@/^1.2pc/[l]\ar@/_1.2pc/[l]
&\D{A_{22}}\ar@/^1.2pc/[l]\ar@/_1.2pc/[l]}\vspace{0.3cm}$$ and
recollements of the rings $A_{ij}$ with $1\le i,j \le 2$, if they
are not derived simple, and so on, until one arrives at derived
simple rings at all positions, or continue to infinitum. All the
derived simple rings appearing in this procedure are called
composition factors of the stratification. The cardinality of the
set of all composition factors (counting the multiplicity) is called
the length of the stratification. If the length of a stratification
is finite, we say that this stratification is finite or of finite
length.

\medskip
{\bf Proof of Corollary \ref{cor1.2}.}  Under the assumption that
$k$ is an algebraically closed field, the following two facts are
known: $(a)$ For any simple regular $R$-module $U$, the algebras
$\End_R(U)$ and  $\End_R(U[\infty])$ are  isomorphic to $k$ and
$k[[x]]$ (see Lemma \ref{lem3.3}(5)), respectively,  and $(b)$ each
tame hereditary algebra having two isomorphism classes of simple
modules is Morita equivalent to the Kronecker algebra.

One the one hand, it follows from Theorem \ref{th1.1} that $\D{B}$
is stratified by $\D{R}$ and $\D{\mathbb{A}_{I}}$, where $I=\{1,2,
\cdots, s\}$ is an index set of the cliques contained in $\mathcal
U$, and the ring $\mathbb{A}_{I}$ is defined in Introduction. Since
$\mathcal{U}$ is a union of finitely many cliques of ${\mathscr S}$,
we know that $\mathbb{A}_{I}$ is equal to $ k((x))^{s}$, the direct
product of $s$ copies of $k((x))$. Thus $\D{\mathbb{A}_{I}}$ has a
stratification by derived module categories with $s$ composition
factors $k((x))$. Note that $\D{R}$ has a stratification by derived
module categories with $r$ copies of the composition factor $k$,
where $r$ is the number of non-isomorphic simple $R$-modules. Thus
$\D{B}$ has a stratification of length $r+ s$ with the composition
factor $k$ of multiplicity $r$, and the composition factor $k((x))$
of multiplicity $s$.

On the other hand, by Corollary \ref{cor4.3}, we know that $\D{B}$
can be stratified by $\D{R_{{\mathcal W}}}$,
$\D{\End_{R_\mathcal{U}}(R_{{\mathcal W}}/R_\mathcal{U})}$ and
$\D{\End_R(R_\mathcal{U}/R)}$, where ${\mathcal W}$ is defined in
Corollary \ref{cor4.3}. Here, we have used the known fact that every
$2\times 2$ triangular matrix ring yields a recollement of derived
module categories of the rings in the diagonal. In the following, we
shall calculate composition factors of $\D{B}$.

First, it follows from Corollary \ref{cor4.3}(3) and Lemma
\ref{lem4.1} that $R_{{\mathcal W}}$ is Morita equivalent to a
Dedekind integral domain and that $\End_{R_\mathcal{U}}(R_{{\mathcal
W}}/R_\mathcal{U})$ is Morita equivalent to $\prod_{j\in
J}T_{c(\mathcal{C}_j)-1}(k)$. It is known from \cite{HKL} that every
Dedekind domain is derived simple. Thus $R_{{\mathcal W}}$
contributes one composition factor to $\D{B}$. It is easy to see
that $\D{T_{c(\mathcal{C}_j)-1}(k)}$ has a stratification with
$c(\mathcal{C}_j)-1$ copies of the composition factor $k$. Thus
$\D{\End_{R_\mathcal{U}}(R_{{\mathcal W}}/R_\mathcal{U})}$ admits  a
stratification with $\sum_{j\in J}\big(c(\mathcal{C}_j)-1\big)$
copies of the composition factor $k$.

Second, combining Corollary \ref{cor4.3}(2) with Corollary
\ref{cor3.4}, we see that $\End_R(R_\mathcal{U}/R)$ is Morita
equivalent to $\prod_{i=1}^s\Gamma\big(c(\mathcal{C}_i)\big)$, where
$\mathcal{U}$ is assumed to be a union of $s$ cliques
$\mathcal{C}_i$ with $1\le i\le s$, and where $\Gamma(m)$ is defined
in Corollary \ref{cor3.4} for each positive integer $m$. Note that
the canonical inclusion $f$ of $\Gamma(m)$ into $M_m(k[[x]])$ is a
ring epimorphism and that $M_m(k[[x]])$ is projective as a left
$\Gamma(m)$-module. Thus the sequence $$0\ra \Gamma(m)\lraf{f}
M_m(k[[x]])\ra \mbox{coker}(f)\ra 0$$ is an
$\add\big(\Gamma(m)E_{m,m}\big)$-split sequence in the category of
all left $\Gamma(m)$-modules, and therefore  $\End_{\Gamma
(m)}\big(\Gamma(m)\oplus M_m(k[[x]])\big)$ and
$\End_{\Gamma(m)}\big(M_m(k[[x]])\oplus \mbox{coker}(f)\big)$ are
derived-equivalent by \cite[Theorem 1.1]{hx2}. Clearly, the former
ring is Morita equivalent to $\Gamma(m)$ and the latter is Morita
equivalent to $\End_{\Gamma(m)}\big(M_m(k[[x]])E_{m,m}\oplus
\mbox{coker}(f)\big)$. Hence $\Gamma(m)$ is derived-equivalent to
$\End_{\Gamma(m)}\big(M_m(k[[x]])\oplus \mbox{coker}(f)\big)$ which
is just the following matrix ring:
$$
{\begin{pmatrix}
k[[x]] &  0  & \cdots &     0\\
k    &  k  & \ddots   &     \vdots\\
\vdots   &\ddots &\ddots & 0\\
k    &  \cdots      & k    & k\\
\end{pmatrix}.}_{m\times m}
$$
For a general consideration of derived equivalences between subrings
of matrix rings, we refer to \cite{ChenYP}.  Thus, we see that
$\D{\Gamma(m)}$ has a stratification with the composition factor
$k[[x]]$ of multiplicity one, and the composition factor $k$ of
multiplicity $m-1$. Therefore, $\D{\End_R(R_\mathcal{U}/R)}$ admits
a stratification with the composition factors: $s$ copies of
$k[[x]]$ and $\sum_{i=1}^s\big(c({\mathcal{C}_i})-1\big)$ copies of
$k$.

Finally, by summarizing up the above discussions, we conclude that
$\D{B}$ has a stratification of length $r+s-1$ with the following
composition factors: $r-2$ copies of $k$, $s$ copies of $k[[x]]$ and
one copy of a fixed Dedekind domain. Here, we use the well known
fact: $\sum_{\mathcal C}\big(c(\mathcal{C})-1\big)=r-2$, where
$\mathcal{C}$ runs over all of the cliques of $R$.  Thus the proof
is completed. $\square$

\bigskip
Let us end this section by mentioning the following questions
suggested by our results.

\smallskip
(1) For tilting modules of the form $R_{\mathcal U}\oplus
R_{\mathcal U}/R$, we have provided a recollement of the derived
categories of their endomorphism rings. It would be interesting to
have a similar result for tilting modules of other types described
in \cite{HJ2}.

(2) In Corollary \ref{cor1.2},  it would be nice to know that
$\D{B}$ has no other composition factors (up to derived equivalence)
except the ones displayed there.

(3) It would be interesting to generalize the results in this paper
to hereditary orders.

(4) Suppose the derived category $\D{A}$ of a ring $A$ admits a
stratification of finite length by derived categories of rings. Does
$\D{A}$ then have only finitely many derived composition factors?
(up to derived equivalence).

\bigskip
 January 6, 2011; Revised: June 20, 2011

\medskip
{\footnotesize
\begin{thebibliography}{99}

\bibitem{HKL}{{\sc L. Angeleri H\"ugel, S. K\"onig} and {\sc Q. Liu}, Recollements
and tilting objects, {\it J. Pure Appl. Algebra} \textbf{215} (2011)
420-438.}

\bibitem{HJ}{{\sc L. Angeleri H\"ugel} and {\sc J. S\'{a}nchez},
Tilting modules arising from ring epimorphisms, {\it Algebr.
Represent. Theory}  \textbf{14} (2011) 217-246. }

\bibitem{HJ2}{{\sc L. Angeleri H\"ugel} and {\sc J. S\'{a}nchez},
Tilting modules over tame hereditary algebras, arXiv:1007.4233v1,
2010.}

%\bibitem{Aus}{{\sc M. Auslander},
%Representation dimension of Artin algebras. Queen Mary College
%Mathematics Notes, Queen Mary College, London, 1971.}

\bibitem{BG}{{\sc K. Bangartz} and {\sc P. Gabriel,} Covering spaces in representation theory, \emph{Invent. Math.} \textbf{65} (1981/82), no.3,
331-378.}

\bibitem{Bz}{{\sc S. Bazzoni}, Equivalences induced by
infinitely generated tilting modules, {\it Proc. Amer. Math. Soc.}
{\bf 138} (2010) 533-544.}

\bibitem{BBD}{{\sc A. A. Beilinson, J. Bernstein} and {\sc P.
Deligne},  Faisceaux pervers, \emph{Asterisque} \textbf{100}
(1982).}

\bibitem{BB} {{\sc S. Brenner} and {\sc M. R. Butler},
Generalizations of the Bernstein-Gelfand-Ponomarev reflection
functors, In: \emph{Representation theory} II. (Eds: V. Dlab and P.
Gabriel), Lecture Notes in Math. \textbf{832}, 103-169, 1980.}

\bibitem{CX}{{\sc H. X. Chen} and {\sc C. C. Xi}, Good tilting
modules and recollements of derived module categories,
arXiv:1012.2176v1, 2010.}

\bibitem{ChenYP}{{\sc Y. P. Chen},  Constructions of derived equivalences, Ph. D. Dissertation, 2011.}

\bibitem{CPS}{{\sc E. Cline, B. Parshall} and {\sc L. Scott},
Algebraic stratification in representation categories, \emph{J.
Algebra} \textbf{117} (1988) 504-521.}

\bibitem{Cohen}{{\sc I. S. Cohen},
On the structure and ideal theory of complete local rings,
\emph{Trans. Amer. Math. Soc.} \textbf{59} (1946) 54-106.}


\bibitem{cohenbook1}{{\sc P. M. Cohn},  \emph{Free rings and their relations},
London Mathematical Society Monographs, \textbf{2}, Academic Press,
London-New York, 1971.}

\bibitem{CB}{{\sc W. W. Crawley-Boevey}, Regular modules for tame hereditary
algebras, \emph{Proc. London Math. Soc.} \textbf{62} (3) (1991)
490-508.}

\bibitem{ct}{{\sc R. Colpi} and {\sc J. Trlifaj}, Tilting modules and tilting torsion
theories, \emph{J. Algebra} \textbf{178} (1995) 614-634.}

%\bibitem{dlab}{{\sc V. Dlab}, The regular representations of the tame
%hereditary algebras, In: S\'eminaire d'Alg\`ebre Paul Dubreil et
%Marie-Paule Malliavin, Lecture Notes in Math. \textbf{1029},
%120-133, 1983.}

\bibitem{dr}{{\sc V. Dlab} and {\sc C. M. Ringel}, Indecomposable representations of graphs and algebras,
{\it Mem. Amer. Math. Soc}. {\bf 173} (1976).}


\bibitem{ep} {{\sc A. J. Engler} and {\sc A. Prestel}, {\it  Valued Fields},
Springer Monographs in Mathematics, Springer-Verlag, Berlin, 2005.}

\bibitem{gt}{{\sc R. G\"obel} and {\sc J .Trlifaj},
{\it  Approximations and endomorphism algebras of modules}, De
Gruyter Expositions in Math. \textbf{41}, Berlin, 2006.}

\bibitem{Happelb}{{\sc D. Happel}, \emph{Triangulated categories in the representation theory of
finite dimensional algebras}, London Math. Soc. Lecture Note Series
\textbf{119}, Cambridge University Press, Cambridge, 1988.}

\bibitem{hr}{{\sc D. Happel} and {\sc C. M. Ringel}, Tilted
algebras, \emph{Trans. Amer. Math. Soc.} \textbf{274} (1982)
399-443.}

\bibitem{hx2}{{\sc W. Hu } and {\sc C. C. Xi},
$\cal D$-split sequences and derived equivalences, {\it Adv. Math.}
\textbf{227} (2011) 292-318.}

%\bibitem{Keller}{{\sc B. Keller}, Derived categories and tilting,
%In: {\it Handbook of tilting theory}, London Mathematical Society
%Lecture Note Series {\bf 332}(2007) 49-104.}

\bibitem{HO}{{\sc H. Krause} and {\sc J. \v S\v tov\'i\v cek}, The Telecope conjecture for hereditary rings
via Ext-orthogonal pairs, {\it Adv. Math.} \textbf {225} (2010)
2341-2364.}

\bibitem{KT}{{\sc P. A. Krylov} and {\sc A. A. Tuganbaev},
\emph{Modules over discrete valuation domains}, De Gruyter
Expositions in Math. \textbf{43}, 2008.}

\bibitem{Lam}{{\sc T. Y. Lam}, \emph{Lectures on  modules and rings}, Graduate Texts in
Math. \textbf{189}, Springer-Verlag, New York, 1999.}

\bibitem{ne} {{\sc J. Neukirch}, {\it  Algebraic number theory},
Grundlehren der MathematischenWissenschaften \textbf{322},
Springer-Verlag, Berlin, 1999.}

\bibitem{RR} {{\sc I. Reiten} and {\sc C. M. Ringel}, Infinite dimensional
representations of canonical algebras, {\it Canad. J. Math.} {\bf
58} (1) (2006) 180-224.}

\bibitem{R0} {{\sc C. M. Ringel}, {\it Representations of $K$-species and bimodules},
\emph{J. Algebra} \textbf{41} (1976) 269-302.}

\bibitem{R} {{\sc C. M. Ringel},
Infinite-dimensional representations of finite-dimensional
hereditary algebras, In: Symposia Mathematica, Vol. XXIII (Conf.
Abelian Groups and their Relationship to the Theory of Modules,
INDAM, Rome, 1977), 321-412, Academic Press, London, 1979.}

\bibitem{R1} {{\sc C. M. Ringel}, {\it Tame algebras and integral
quadratic forms}, Lecture Notes in Math. {\bf 1099},
Springer-Verlag, Berlin,  1984.}

\bibitem{Sch1}{{\sc{A. Schofield}}, {\it Representations of rings over
skew fields}, London Math. Soc. Lecture Note Series {\bf 2},
Cambridge University Press, Cambridge, 1985.}

\bibitem{Sch2}{{\sc{A. Schofield}}, Universal localization
for hereditary rings and quivers, In: \emph{Ring Theory} (Antwerp,
1985), Lecture Notes in Math. {\bf 1197}, 149-164, 1986.}

\end{thebibliography}}
\end{document}